\documentclass[letterpaper,11pt,leqno,color]{amsart}
\usepackage[latin1]{inputenc}
\usepackage{latexsym}
\usepackage[dvips]{graphicx}
\usepackage{amsfonts}
\usepackage{epsfig}
\usepackage{geometry}
\usepackage{color}
\usepackage{amsmath}
\usepackage{amssymb}

\newtheorem{thm}{Theorem}
\newtheorem{prop}[thm]{Proposition}
\newtheorem{lem}[thm]{Lemma}
\newtheorem{cor}[thm]{Corollary}
\newtheorem{claim}[thm]{Claim}
\newtheorem{rem}[thm]{Remark}

\newcommand{\N}{\mathbf{N}}

\newcommand{\R}{\mathbf{R}}

\newcommand{\Rcupinf}{\R\cup\{+\infty\}}

\newcommand{\luno}{\ell^1}
\newcommand{\ltwo}{\ell^2}

\renewcommand{\H}{\mathcal{H}}

\newcommand{\M}{\mathcal{M}}

\renewcommand{\S}{\mathcal{S}}

\newcommand{\eps}{\varepsilon}

\newcommand{\limty}[1]{\lim\limits_{#1\to\infty}}
\newcommand{\limn}{\limty{n}}

\newcommand{\dem}{\bigskip\noindent{\bf Proof. }}
\newcommand{\bx}{\hfill$\blacksquare$\\ \bigskip}

\def\argmin{\hbox{Argmin}}

\def\dom{\hbox{dom}}

\newcommand{\cali}{\mathcal}

\def\noi{\noindent}
\def\u{\underline}

\begin{document}

\title{Splitting forward-backward penalty scheme for constrained variational problems}

\author{Marc-Olivier Czarnecki}

\author{Nahla Noun}

\author{Juan Peypouquet}

\address{Institut de Math\'ematiques et Mod\'elisation de Montpellier, UMR 5149 CNRS, Universit\'e Montpellier 2, place Eug\`ene Bataillon,
34095 Montpellier cedex 5, France}
\email{marco@univ-montp2.fr, nahla.noun@yahoo.fr}

\address{Departamento de Matem\'atica, Universidad T\'ecnica Federico Santa Mar\'\i a, Avenida  Espa\~na 1680, Valpara\'\i so, Chile.}
\email{juan.peypouquet@usm.cl}
\date{\today}

\maketitle

\begin{abstract}

We study a forward backward splitting algorithm that solves the variational inequality
\begin{equation*}
A x +\nabla \Phi(x)+ N_C (x) \ni 0
\end{equation*}
where $\H$ is a real Hilbert space,  $A: \H\rightrightarrows \H$ is a maximal monotone operator, $\Phi: \H\to\R$ is a smooth convex function, and $N_C$ is the outward normal cone to a closed convex set $C\subset\H$. The constraint set $C$ is represented as the intersection of the sets of minima of two convex penalization function $\Psi_1:\H\to\R$ and $\Psi_2:\H\to\R\cup \{+\infty\}$. The function $\Psi_1$ is smooth, the function $\Psi_2$ is proper and lower semicontinuous.
Given a sequence $(\beta_n)$ of  penalization parameters which tends to infinity, and a sequence of positive time steps $(\lambda_n)$, the algorithm
\begin{equation*}
\ \left\{\begin{array}{rcl}
x_1 & \in & \H,\\
x_{n+1} & = & (I+\lambda_n A+\lambda_n\beta_n\partial\Psi_2)^{-1}(x_n-\lambda_n\nabla\Phi(x_n)-\lambda_n\beta_n\nabla\Psi_1(x_n)),\  n\geq 1.
\end{array}\right.\leqno(SFBP)
\end{equation*}
performs forward steps on the smooth parts and backward steps on the other parts. 
Under suitable assumptions, we obtain weak ergodic convergence of the sequence $(x_n)$ to a solution of the variational inequality. Convergence is strong when either $A$ is strongly monotone or $\Phi$ is strongly convex.  We also obtain weak convergence of the whole sequence $(x_n)$ when $A$ is the subdifferential of a proper lower-semicontinuous convex function.
This provides a unified setting for several classical and more recent results, in the line of historical research on continuous and discrete gradient-like systems. 
\end{abstract}


\vspace{1cm}


%
%
%
%
%
%
%
%
%
%
%

\section{Introduction}

In 1974, R. Bruck  \cite{Bruck} showed that the trajectories of the steepest descent system
\begin{equation*}
\dot{x}(t)+\partial \Phi(x(t))\ni 0
\end{equation*}
minimize the convex, proper, lower-semicontinuous potential $\Phi$ defined on  a real Hilbert space $\H$. They weakly converge to a point in the minima of  $\Phi$ and the potential decreases along the trajectory toward its minimal value, provided $\Phi$ attains its minimum. 
When the semigroup is generated by the differential inclusion
\begin{equation*}
\dot{x}(t)+A(x(t))\ni 0
\end{equation*}
with a maximal monotone operator $A$ from  $\H$ to $\H$, J.-B. Baillon and H. Br\'ezis \cite{BaiBre} provided in 1976 the convergence in average to an equilibrium of $A$.
These results are sharp. If the operator is a rotation in $\R^2$, the trajectories do not converge, except the stationary one. 
J.-B. Baillon \cite{Bai}, provided an example in 1978, where the trajectories of the steepest descent system do not strongly converge, although they weakly converge. In some sense, his example is an extension to a Hilbert space of the  rotation in $\R^2$.
The keytool for the proof of these two results is Opial's lemma \cite{Opi} that gives  weak convergence without a priori knowledge of the limit.
In 1996, H. Attouch and R. Cominetti \cite{AttCom} coupled approximation methods with the steepest descent system, in  particular by adding a Tikhonov regularizing term:
\begin{equation*}
\dot{x}(t)+\partial \Psi(x(t))+\varepsilon (t) x(t)\ni 0.
\end{equation*}
The parameter $\varepsilon$ tends to zero and the potential $\Psi$ satisfies usual assumptions. As it yields the steepest descent system for  $\varepsilon=0$, one can expect the trajectories to weakly converge.
The striking point of their results is the strong convergence of the trajectories, when $\varepsilon$ tends to $0$ slowly enough, that is $\varepsilon$ does not belong to $l^1$. Then the strong limit is the point of minimal norm among the minima of $\Psi$.
This seems rather surprising at first: without a regularizing term $\varepsilon (t) x(t)$, we know that we only have weak convergence, and with the regularizing term, convergence is strong.
We propose the following explanation: set 
$
\Phi(x)= \frac{1}{2}\|x\|^2
$
so that the regularizing term  writes $\varepsilon (t) x(t)= \varepsilon (t)\nabla \Phi(x(t))$. Then, by a change of time, valid for $\varepsilon\notin l ^1$, see~\cite{AC}, we can reformulate the system as a penalized system
\begin{equation*}
\dot{x}(t)+\nabla \Phi(x(t))+\beta(t)\partial \Psi(x(t))\ni 0,
\end{equation*}
with a parameter $\beta$ that tends to infinity.
But now we are looking at a steepest descent system for the stongly convex function $\Phi$ with a penalization potential $\Psi$, possibly equal to $0$. And it is known that the trajectories of the steepest descent system strongly converge when the potential is strongly convex. In 2004, A. Cabot \cite{Cab1} generalized part of Attouch and Cominetti's result to the case of a strongly convex potential. The penalization acts as a constraint and forces the limit to belong to the minima of $\Psi$. 
It appeared natural to add a penalization, rather than a perturbation or a regularization, to the first order differential inclusion with a maximal monotone operator, and to the  steepest descent system with a -not necessarily strongly- convex potential. Moreover, penalization methods enjoy gret practical interest when others, such as projection methods, encounter intrinsic implementation difficulties, for example when the constraint set is given by nonlinear inequalities.
H. Attouch and M.-O. Czarnecki \cite{AC} showed in 2010 that the trajectories of
\begin{equation*}
\dot{x}(t) +  A(x(t)) +  \beta (t) \partial \Psi (x(t))\ni 0.
\end{equation*}
weakly converge in average to a constrained equilibrium
\begin{equation*}
x_{\infty}\in(A + N_C)^{-1}(0),
\end{equation*}
that convergence is strong when $A$ is strongly monotone, and in the subdifferential case, the trajectories of
\begin{equation*}
\dot{x}(t) +  \partial \Phi(x(t)) +  \beta (t) \partial \Psi (x(t))\ni 0.
\end{equation*}
weakly converge to a constrained minimum
\begin{equation*}
x_{\infty}\in\mbox{argmin}\{ \Phi|\mbox{argmin}\Psi\}.
\end{equation*}
Besides assuming the parameter $\beta$ to tend to $+\infty$, their main assumption relates the geometry of the penalization potential $\Psi$ to the growth rate of the penalization parameter $\beta$, namely
\begin{equation*}
 \int_{0}^{+\infty} \beta (t) \left[\Psi^*
  \left(\frac{p}{ \beta (t)}\right) - \sigma_C \left(\frac{p}{ \beta
    (t)}\right)\right]dt < + \infty
\end{equation*}
for every $p$ in the range of $N_C$, the normal cone to $C=\mbox{argmin}\Psi$.  Here $\Psi^*$ denotes the Fenchel conjugate of $\Psi$ and $ \sigma_C$ the support function of $C$. A detailed analysis of the condition is done in~\cite{AC}. Let us just mention that, when $\Psi =  \frac{1}{2} {\rm dist}_C^2$, it reduces to
$\int_{0}^{+\infty} \frac {1}{\beta (t)} dt < + \infty$.
When $\Psi=0$, then $C=\H$ and the only $p$ in the normal cone is $0$, and the condition is fulfilled. So one recovers the results of Bruck and of Baillon and Br\'ezis.

\subsection*{Discretization} In order to compute the trajectories of the system, and obtain workable algorithms, we need to discretize it.
The implicit discretization of unpenalized first order system is the famous proximal algorithm, well studied from the mid seventies (B. Martinet~\cite{Mar1} and~\cite{Mar2}, R.T. Rockafellar~\cite{Rockafellar}, H. Br\'ezis and P.L. Lions.~\cite{bresis-lions},\dots):
\begin{equation*}
x_{n+1}=(I+\lambda_n A)^{-1}x_n.
\end{equation*}
Following the same path,  in 2011, Attouch, Czarnecki and Peypouquet \cite{algo1} discretized the penalized continuous system implicitly to obtain the backward algorithm:
\begin{equation*}
x_{n+1}=(I+\lambda_n A+\lambda_n\beta_n\partial\Psi)^{-1}x_n.
\end{equation*}
They provide weak convergence in average to a constrained equilibrium
$
x_{\infty}\in(A + N_C)^{-1}(0),
$
 strong convergence when $A$ is strongly monotone, and weak convergence in the subdifferential case. A basic assumption on the time step is $(\lambda_n)\notin l^1$, which corresponds to $t\to+\infty$ in continuous time. The key assumption is the discrete counterpart of the assumption in continuous time:
\begin{equation*}
\sum\limits_{n=1}^\infty
\lambda_n\beta_n\left[\Psi^*\left(\dfrac{z}{\beta_n}\right)-\sigma_C\left(\dfrac{z}{\beta_n}\right)\right]<\infty
\end{equation*}
for every $p$ in the range of $N_C$. Again, in the case where $\Psi=0$, one recovers classical results for the proximal algorithm. The main drawback of the backward algorithm is the cost of every step in the computation of the discrete trajectory. The explicit discretization of the steepest descent system can be traced back to A. Cauchy~\cite{Cau} in 1847, who gave indeed the idea of the discrete gradient method
\begin{equation*}
x_{n+1}=x_n-\lambda_n\nabla\Phi(x_n),
\end{equation*}
with no proof of convergence, and before the continuous steepest descent system. J. Peypouquet~\cite{algo4} discretized the continuous penalized system explicitly in 2012 to obtain a forward algorithm, in a regular setting
\begin{equation*}
x_{n+1}=x_n-\lambda_n\nabla\Phi(x_n)-\lambda_n\beta_n\nabla\Psi(x_n).
\end{equation*}
He shows weak convergence of the trajectories to a constrained minimum
\begin{equation*}
x_{\infty}\in\mbox{argmin}\{ \Phi|\mbox{argmin}\Psi\}
\end{equation*}
provided the gradient of the potential $\Phi$ is Lipschitz continuous. Together with the key assumption described before relating the Fenchel conjugate of $\Psi$ and the sequences $(\lambda_n)$ and $(\beta_n)$, he requires an assumption combining a bound on these sequences and the Lipschitz constant of $\nabla\Psi$. It slightly differs, and is a consequence of
\begin{equation*}
\limsup \limits_{n\to\infty}\lambda_n\beta_n<\frac{2}{ L_{\nabla \Psi}}.
\end{equation*}
To make things short, forward algorithms are more performing, but require more regularity assumptions and convergence is more complicated to prove, while backward algorithms apply to more general cases, convergence is easier to prove
 but they are not so efficient.
As the constraint set can possibly be described by a regular penalization function, the next idea, developed in~\cite{algo3}, is to perform a forward step on the regular part $\Psi$, and a backward step on the other part to obtain the forward-backward algorithm:
\begin{equation*}
x_{n+1}=(I+\lambda_n A)^{-1}(x_n-\lambda_n\beta_n\nabla\Psi(x_n)).
\end{equation*}
We obtain again weak convergence in average to a constrained equilibrium
$
x_{\infty}\in(A + N_C)^{-1}(0),
$
 strong convergence when $A$ is strongly monotone, and weak convergence in the subdifferential case. Together with the key summability assumption relating $\Psi$ and the parameters $(\lambda_n)$ and $(\beta_n)$, we assume regularity of the function $\Psi$, that it is differentiable with a Lipschitz gradient.
We need the same assumption on $(\lambda_n)$, $(\beta_n)$, and the Lipschitz constant of $\nabla\Psi$:
\begin{equation*}
\limsup \limits_{n\to\infty}\lambda_n\beta_n<\frac{2}{ L_{\nabla \Psi}}.
\end{equation*}
The bound is strict in general, and is close to being sharp: equality, with a precise potential $\Psi$, corresponds to a result of Passty~\cite{Pas_1979} on alternate algorithms. A detailed analysis is given in~\cite{algo3}.

\subsection*{Regularity based splitting} We now have three different algorithms depending on the regularity of the data: backward algorithms, forward algorithms, forward-backward algorithms. Convergence holds under the same summability assumption involving the Fenchel conjugate of $\Psi$, and similar regularity assumptions to perform a forward step.

What if the maximal monotone operator has a regular part, and if the penalization potential decomposes with a regular part? Can we guarantee convergence if we perform forward steps on the regular parts, while keeping the backward step on the other parts? Can we provide a unified setting for the previous algorithms?

Let $A$ be a maximal monotone operator, let $\Phi$, $\Psi_1$, $\Psi_2$ be convex proper lower semicontinuous potentials.  The functions $\Phi$ and $\Psi_2$ are defined everywhere and differentiable with a Lipschitz gradient. Set $C=\argmin\Psi_1\cap\argmin\Psi_2$, which corresponds to the decomposition of the penalization function as the sum of a smooth part and a general part, and assume that $C$ is not empty.
 The penalized system
\begin{equation}\label{mag}
\dot{x}(t) +  (A+\nabla \Phi)(x(t)) +  \beta (t) (\partial \Psi_1+\nabla \Psi_2)  (x(t))\ni 0.
\end{equation}
allows to solve the variational inequality
\begin{equation}\label{ineq}
 0\in Ax+ \nabla \Phi(x)+N_C(x).
 \end{equation}
We discretize this last continuous penalized system by making a forward step on the regular parts $\Phi$ and $\Psi_2$, and a backward step on $A$ and  $\Psi_1$.
Given a positive sequence $(\beta_n)$ of {\em penalization parameters}, along with a positive sequence $(\lambda_n)$ of {\em step sizes}, we consider the \emph{splitting forward-backward penalty algorithm} (SFBP), defined as follows:
\begin{equation*}
\ \left\{\begin{array}{rcl}
x_1 & \in & \H,\\
x_{n+1} & = & (I+\lambda_n A+\lambda_n\beta_n\partial\Psi_2)^{-1}(x_n-\lambda_n\nabla\Phi(x_n)-\lambda_n\beta_n\nabla\Psi_1(x_n)),\  n\geq 1.
\end{array}\right.\leqno(SFBP)
\end{equation*}

Under suitable  assumptions, including the expected geometrical condition involving the Fenchel conjugate of the penalization potential and the expected relationship between the parameters and the Lipschitz constant  of $\nabla \Psi_1$, we prove that, as $n\to\infty$, the sequences generated by the (SFBP) algorithm converge to a constrained equilibrium in $S=(A+ \nabla \Phi+N_C)^{-1}0$
\begin{itemize}
    \item [i)] weakly in average if $A$ is any maximal monotone operator (Theorem \ref{theo2'});
    \item [ii)] strongly if $A$ is strongly monotone, or if $\Phi$ is strongly convex (Theorem \ref{theo1'});
    \item [iii)] weakly if $A$ is the subdifferential of a proper, lower-semicontinuous and convex function (Theorem \ref{theo2}).
\end{itemize}

Besides its applicability to problems that combine smooth and nonsmooth features, the (SFBP) algorithm allows us to study, in a unified framework, the classical and more recent methods to solve constrained variational inequalities.

If $\Phi=\Psi_1=\Psi_2\equiv 0$, we recover the {\em proximal point algorithm} . If $\Phi=\Psi_1\equiv 0$, the (SFBP) algorithm corresponds to the purely implicit {\em prox-penalization algorithm} from \cite{algo1}.
 The {\em gradient method},  is recovered in the case where $A\equiv 0$ and $\Psi_1=\Psi_2\equiv 0$. 
 If $A\equiv 0$ and $\Psi_2\equiv 0$, we obtain the purely explicit {\em diagonal gradient scheme} from \cite{algo4}.
We get the {\em forward-backward splitting} (a combination of the {\em proximal point algorithm} and of the {\em gradient method}, see \cite{Pas_1979}) if $\Psi_1=\Psi_2\equiv 0$. The case $\Phi=\Psi_2\equiv 0$ gives a semi-implicit penalty splitting method studied in \cite{algo3} and \cite{Noun-Pey}. 

\subsection*{Organization of the paper}

The paper is organized as follows: We begin by describing and commenting the hypothesis, and stating the main theoretical results, in Section~\ref{S:statements}. Next, in Section~\ref{S:comparison}, we present several special cases, and compare them with classical and more recent methods to solve constrained optimization problems. In particular, this work extends and unifies some previous developments progressively achieved by our work group.
In Section~\ref{secexamples}, we describe a model for the sparse-optimal control of a linear system of ODE's. It illustrate how the decomposition of objective potential and penalization naturally arises. Finally, we present the proofs in several steps  Section~\ref{secproposition}.
\if{ 
Weak ergodic convergence to a solution of inclusion \eqref{ineq} is proved in Section \ref{S:ergodic}. The strongly convergent case is studied in Section \ref{S:strong}. In Section \ref{S:weak}, we establish the weak convergence of the whole sequence to a point in the set of constrained minima,  when $A$ is the subdifferential of a proper, lower-semicontinuous and convex function. Finally, 
}\fi

\section{Main results}\label{S:statements}\label{secresults}

Let $\H$ be a real Hilbert space. We first recall some facts about convex analysis and maximal monotone operator theory.
Let $\Gamma_0(\H)$ be the set of all proper (not identically equal to $+\infty$) lower-semicontinuous convex functions from $\H$ to $\R\cup \{+\infty\}$.
Given
$F\in\Gamma_0(\H)$ and $x\in\H$, the {\em subdifferential} of $F$ at $x$ is the set
$$\partial F(x)=\{x^*\in \H:F(y)\ge F(x)+\langle x^*,y-x\rangle\hbox{  for all  }y\in\H\}.$$

Given a nonempty closed convex set $C\subset\H,$ its {\em indicator function} is defined as $\delta_C(x)=0$ if $x\in C$ and $+\infty$ otherwise.
The {\em normal cone} to $C$ at $x$ is
$$N_C(x)=\{x^*\in\H:\langle x^*,c-x\rangle\le 0\hbox{  for all  }c\in C\}$$ if $x\in C$ and $\emptyset$ otherwise. Observe that 
$\partial\delta_C=N_C$. 
A \textit{monotone operator} is a set-valued mapping $A:\H\rightarrow\H$ such that $\langle x^*-y^*,x-y\rangle\geq 0$
whenever $x^*\in Ax$ and $y^*\in Ay$.
It is \textit{maximal monotone} if its graph is not properly contained in the graph of any other monotone operator.
It is convenient to identify a maximal monotone operator $A$ with its graph, thus we equivalently write $x^*\in Ax$ or $[x,x^*]\in A$.
The inverse $A^{-1}:\H\rightarrow\H$ of $A$ is defined by $x\in A^{-1}x^*\Leftrightarrow x^*\in Ax$.
It is still a maximal monotone operator.
For any maximal monotone operator $A:\H\rightarrow\H$ and for any $\lambda>0$, the operator $I+\lambda A$ is surjective by Minty's Theorem
(see \cite{brezis} or \cite{Pey_Sor_2010}). The operator $(I+\lambda A)^{-1}$ is nonexpansive and everywhere defined.
It is called the \textit{resolvent} of $A$ of index $\lambda$.

Finally recall that the subdifferential of a function in $\Gamma_0(\H)$ is maximal monotone.

\subsection{Assumptions}
Let $A$ be a maximal monotone operator, let $\Phi$, $\Psi_1$, $\Psi_2$ be convex proper lower semicontinuous potentials with
$$C=\hbox{Argmin$\Psi_1\ \cap$ Argmin$\Psi_2$} \neq \emptyset.$$
 The functions $\Phi$ and $\Psi_2$ are defined everywhere and differentiable with a Lipschitz gradient.
The {\em Fenchel conjugate} of a proper, lower-semicontinuous and convex function $F:\H\to\Rcupinf$ is the function $F^*:\H\to\Rcupinf$ defined by
$$F^*(x^*)=\sup_{y\in\H}\left\{\langle y,x^*\rangle -F(y)\right\}$$
for each $x^*\in\H$. It is also proper, lower-semicontinuous and convex. Given a nonempty, closed and convex set $C$, its {\em support function} $\sigma_C$ is defined as $\sigma_C(x^*)=\sup\limits_{c\in C}\langle x^*,c\rangle$ for $x^*\in\H$. Observe that $\delta_C^*=\sigma_C$. Notice also that $x^*\in N_C(x)$ if, and only if, $\sigma_C(x^*)=\langle x^*,x\rangle$.\\

The main set of hypotheses is the following:
$$\left\{\begin{array}{rl}
i)&{\mathbf T}=A+\nabla\Phi+N_C\ \hbox{is maximal monotone and}\ S={\mathbf T}^{-1}(0)\neq \emptyset ;\\
ii)& \nabla\Phi\ \hbox{is}\ L_{\Phi}\hbox{-Lipschitz-continuous\ and}\ \nabla\Psi_1 \ \hbox{is}\ L_{\Psi_1}\hbox{-Lipschitz-continuous};\\
iii)&\hbox{For each $z\in N_C(\H)$},\ \ \  \sum\limits_{n=1}^\infty
\lambda_n\beta_n\left[\left(\Psi_1+\Psi_2\right)^*\left(\dfrac{z}{\beta_n}\right)-\sigma_C\left(\dfrac{z}{\beta_n}\right)\right]<\infty;\\
iv)& \sum\limits_{n=1}^\infty\lambda_n=+\infty,\ \ \sum\limits_{n=1}^\infty L_\Phi\lambda^2_n<+\infty,\ \ \hbox{and}\ \
\limsup\limits_{n\to\infty}( L_{\Psi_1}\lambda_n\beta_n)<2.
\end{array}\right.\leqno{\rm (H_0)}$$

We adress some remarks on Hypothesis ${\rm (H_0)}$ in order.

\subsubsection*{On Part i)}
It is a well-posedness and qualification condition ensuring that
\begin{equation}
\label{E:charact}
\bar u\in S\quad\hbox{if, and only if,}\quad \langle w,u-\bar u\rangle \ge 0\hbox{ for all }[u,w]\in {\mathbf T}.
\end{equation}
	If $A$ is the subdifferential of a proper, lower-semicontinuous and convex function $\Phi_2 $, the maximal monotonicity of  ${\mathbf T}$ implies
\begin{equation}\label{E:S_subdifferential}
S=\hbox{Argmin}\{\Phi(x)+\Phi_2(x): x\in C\}.
\end{equation}
In this situation, $S$ can be interpreted as the set of solutions of a hierarchical optimization problem, where $\Phi+\Phi_2$ and $\Psi_1+\Psi_2$ are
primary and secondary criteria, respectively.
In this case, maximality of ${\mathbf T}$ holds under some qualification condition, such as Moreau-Rockafellar or Attouch-Br\'ezis.
\subsubsection*{On Part ii)}
It is standard for the convergence of gradient-related methods (see \cite{bertsekas}).
\subsubsection*{On Part iii)}
It was introduced in \cite{algo1}, following \cite{AC}.
The potentials $\Psi_1$ and $\Psi_2$ enter the algorithm only via their subdifferentials. Thus it is not a restriction to assume $\min \Psi_1=\min \Psi_2=0$. Otherwise, one should replace $\Psi_i$ by $\Psi_i-\min \Psi_i$ in the corresponding statements.
 In the unconstrained case ($\Psi_1=\Psi_2\equiv 0$), the condition is trivially satisfied since $N_C(\H)=\{0\}$, $\Psi^*(0)=0$ and $\sigma_{\H}(0)=0$. We refer to  \cite{algo1} for discussion and
 sufficient conditions.
Note that the constraint set $C$ is the set of minima of the potential $\Psi_1+\Psi_2$, which leads naturally to an assumption on the Fenchel conjugate of the sum $\Psi_1+\Psi_2$, and involving points in the normal cone $N_C$. In our setting, considering alternatively the two separate corresponding conditions on $\Psi_1^*$ and $\Psi_2^*$
would require extra qualification conditions.

\subsubsection*{On Part iv)}
The nonsummability condition in Part iv) is standard for the proximal point algorithm (see \cite{bresis-lions}) and gradient-related methods (see \cite{bertsekas}). The second condition holds if either $\Phi$ is affine (that is $L_\Phi=0$) or $(\lambda_n)$ is in $\ltwo$. We write $\limsup\limits_{n\to\infty}(L_{\Psi_1}\lambda_n\beta_n)<2$ rather than $\limsup\limits_{n\to\infty}\lambda_n\beta_n<\frac{2}{ L_{\Psi_1}}$ to include the case where  $L_{\Psi_1}=0$ ($\Psi_1\equiv 0$).\\

\subsection{Convergence results}
Take a sequence $(x_n)$ generated by 
the \emph{splitting forward-backward penalty algorithm} (SFBP):
\begin{equation*}
\ \left\{\begin{array}{rcl}
x_1 & \in & \H,\\
x_{n+1} & = & (I+\lambda_n A+\lambda_n\beta_n\partial\Psi_2)^{-1}(x_n-\lambda_n\nabla\Phi(x_n)-\lambda_n\beta_n\nabla\Psi_1(x_n)),\  n\geq 1,
\end{array}\right.\leqno(SFBP)
\end{equation*}
which corresponds to the implicit-explicit discretization 
$$
\frac{x_n-x_{n+1}}{\lambda_n}-\nabla\Phi(x_n)-\beta_n\nabla\Psi_1(x_n) \in Ax_{n+1}+\beta_n\partial\Psi_2(x_{n+1}).
$$
of the penalized differential inclusion \eqref{mag}.
We do not discuss in detail the existence of trajectories. Maximality of $A+\beta_n\partial \Psi_2$ for all $n\in \N$ is sufficient in view of Minty's theorem.  Moreover, according the discussion in Subsection \ref{SS:inexact}, it is possible to consider the above inclusion replacing the subdifferential operator by some enlargement, such as the {\em $\eps$-approximate subdifferential}.\\

The kind of convergence depends on the nature of the operator $A$.\\

When $A$ is any maximal monotone operator we prove the weak ergodic convergence of the algorithm to a point in $S$. More precisely, let $(x_n)$ be a sequence generated by (SFBP) and let $\tau_n=\sum\limits_{k=1}^{n}\lambda_k$.
We define the following sequences of weighted averages:
\begin{equation}\label{zn}
z_n=\dfrac{1}{\tau_n}\sum\limits_{k=1}^{n}\lambda_kx_k~\hspace*{3cm}\widehat{z}_n=\dfrac{1}{\tau_n}\sum\limits_{k=1}^{n}\lambda_kx_{k+1}.
\end{equation}
Although quite similar, they converge under slightly different assumptions, as we shall see.

\begin{thm}\label{theo2'}
Assume that ${\rm (H_0)}$ holds. Then the sequence $(\widehat{z}_n)$ converges weakly to a point in $S$ as $n\to\infty$. If we moreover assume that $(\lambda_n)\in\ell^2$ then the sequence $(z_n)$
converges weakly to a point in $S$ as $n\to\infty$.
\end{thm}

Under further assumptions, it is possible to obtain strong or weak convergence of {\em the whole sequence} $(x_n)$. Recall that $A$ is is \textit{strongly monotone} with parameter $\alpha>0$ if
$$\langle x^*-y^*,x-y\rangle\geq\alpha\|x-y\|^2 $$
whenever $x^*\in Ax$ and $y^*\in Ay$. The function $\Phi$ is strongly convex if $\nabla \Phi$ is  strongly monotone.
The set of zeros of a maximal monotone operator which is strongly monotone must contain exactly one element. We have the following:

\begin{thm}\label{theo1'}
Let ${\rm (H_0)}$ hold. If the operator $A$ is strongly monotone, or if the potential $\Phi$ is strongly convex, then every sequence $(x_n)$ generated by algorithm (SFBP) converges strongly to the unique $u\in S$ as $n\to\infty$.
\end{thm}

Finally, a function $F:\H\to\Rcupinf$ is {\em boundedly inf-compact} if the sets of the form
$$\{\,x\in\H\ :\ \|x\|\le R,\hbox{ and }F(x)\le M\,\},$$
are relatively compact for every $R\geq 0$ and $M\in \R$.\\

We shall prove that if $A$ is the subdifferential of a proper, lower-semicontinuous and convex function $\Phi_2$, weak convergence of the sequences generated by the (SFBP) algorithm can be guaranteed if either $\Phi_2$ is boundedly inf-compact, the penalization parameters satisfy a {\em subexponential} growth condition, or in the unconstrained case. More precisely, we have the following:

\begin{thm}\label{theo2}
Let ${\rm (H_0)}$ hold with $A=\partial \Phi_2$. Assume that any of the following conditions holds:
\begin{itemize}
\item[(i)] $\liminf\limits_{n\to\infty}\lambda_n\beta_n>0$ and the function $(\Phi+\Phi_2)$ or $(\Psi_1+\Psi_2)$ is boundedly inf-compact;
\item[(ii)] $\liminf\limits_{n\to\infty}\lambda_n\beta_n>0,$ $(\lambda_n)$ is bounded and $\beta_{n+1}-\beta_n\leq K\lambda_{n+1}\beta_{n+1}$
for some $K>0$; or
\item[(iii)] $\Psi_1=\Psi_2=0.$
\end{itemize}
Then, the sequence $(x_n)$ converges weakly to a point in $S$ as $n\to\infty$.
Moreover, convergence is strong in case $(i)$ if $\Psi_1+\Psi_2$ is  boundedly inf-compact.\\

The sequence $(x_n)$  minimizes $\Phi+\Phi_2$ in cases $(ii)$ and $(iii)$:
 $$\lim\limits_{n\to\infty}(\Phi+\Phi_2)(x_n)=\min_{C}(\Phi+\Phi_2).$$
\end{thm}

The proofs of Theorems \ref{theo2'}, \ref{theo1'} and \ref{theo2} will be completed in Sections \ref{S:ergodic}, \ref{S:strong} and \ref{S:weak}, respectively.
One cannot expect to have strong convergence in Theorem \ref{theo2} in general, see the comment after Corollary~\ref{cor140117}.\\

\subsection{Inexact computation of the iterates}\label{SS:inexact}

Convergence also holds if the iterates are computed inexactly provided the errors are small enough. More precisely, consider the {\em inexact splitting forward-backward penalty algorithm} given by

$$\ \left\{\begin{array}{rcl}
x_1 & \in & \H,\\
x_{n+1} & = & (I+\lambda_n A+\lambda_n\beta_n\partial\Psi_2)^{-1}(x_n-\lambda_n\nabla\Phi(x_n)-\lambda_n\beta_n\nabla\Psi_1(x_n)-\zeta_n)+\xi_n,\\
&& n\geq 1.
\end{array}\right.\leqno(SFBP_\eps)$$

We recall the following result from \cite{Alv_Pey_2010}:

\begin{lem}
Let $(P_n)$ be a sequence of nonexpansive functions from $\H$ into $\H$. Let $(\eps_n)$ be a positive sequence in $\luno$. If every sequence $(x_n)$ satisfying
$$x_{n+1}=P_n(x_n),\quad n\ge 1$$
converges weakly (resp. strongly, resp. weakly or strongly in average), then the same is true for every sequence $(\tilde x_n)$ satisfying
$$\|\tilde x_{n+1}-P_n(\tilde x_n)\|\le\eps_n,\quad n\ge 1.$$
\end{lem}

Following the arguments in the proof of \cite[Proposition 6.3]{algo1}, we obtain

\begin{cor}
Let $(\zeta_n)$ and $(\xi_n)$ be nonnegative sequences in $\luno$, and let $(x_n)$ verify ${\rm (SFBP_\eps)}$. Then Theorems \ref{theo2'}, \ref{theo1'} and \ref{theo2} remain true.
\end{cor}

\subsection{Forward Backward Backward algorithm and full splitting}

The (SFBP) algorithm is a step toward full splitting. It allows to understand the different roles played by the regular parts -allowing for forward steps- and the general parts -needing backward steps. It requires to compute the resolvent of the sum of two maximal monotone operators, which may be a hard task. The full splitting af the backward step is achieved in \cite{algo1}. Following the same path, let us define the \emph{splitting forward-backward-backward penalty algorithm} (SFBBP), as follows:
\begin{equation*}
\ \left\{\begin{array}{rcl}
x_1 & \in & \H,\\
x_{n+1} & = & (I+\lambda_n\beta_n\partial\Psi_2)^{-1} (I+\lambda_n A)^{-1}(x_n-\lambda_n\nabla\Phi(x_n)-\lambda_n\beta_n\nabla\Psi_1(x_n)),\  n\geq 1.
\end{array}\right.\leqno(SFBBP)
\end{equation*}
The complete study of (SFBBP) goes beyond the scope of this paper. We believe that the convergence results should hold, by making use of the techniques in 
\cite{algo1}.

\section{Comparison with classical and more recent methods}\label{S:comparison}\label{seccomparison}

In this section we examine some particular cases, where some of the functions or the operator involved in \eqref{ineq} vanish.

\subsection{The backward algorithm}
Taking the two potentials $\Phi$ and $\Psi_1$ to be zero, the forward step disappears and the  (SFBP) algorithm turns into a purely backward algorithm.

\subsubsection{The unconstrained case : the proximal point algorithm  }

If additionnally  $\Psi_2$ is zero, we obtain the {\em proximal point algorithm}:
$$\left\{\begin{array}{rcl}
x_1 & \in & H,\\
x_{n+1} & = & (I+\lambda_n A)^{-1}(x_n)\ \ {\rm for\ all}\ n\geq 1.
\end{array}\right.\leqno{(PROX)}$$

This method was originally introduced in \cite{Mar1}, using the idea of {\em proximity operator} from \cite{Moreau}. It was further developed in \cite{Rockafellar}, \cite{bresis-lions} and \cite{lions}. Its popularity is due to the fact that, despite its iteration-complexity, convergence can be granted under minimal hypotheses.\\

Let $(x_n)$ be a sequence generated by (PROX) and define $(z_n)$ and $(\widehat z_n)$ and in \eqref{zn}.

\begin{cor}\label{cor140117}
Let $S\neq\emptyset$ and $(\lambda_n)\notin\luno$. As $n\to\infty$, we have the following:
\begin{itemize}
	\item[i)] The sequence $(\widehat{z}_n)$ converges weakly to a point in $S$;
	\item[ii)] If $(\lambda_n)\in\ell^2$, then the sequence $(z_n)$ converges weakly to a point in $S$;
	\item[iii)] If $A$ is strongly monotone, then $(x_n)$ converges strongly to the unique point in $S$; and
	\item[iv)] If $A=\partial \Phi_2$, then $(x_n)$ converges weakly to a point in $S$, with $\lim\limits_{n\to\infty}\Phi_2(x_n)=\min(\Phi_2)$.
\end{itemize}
\end{cor}

Part i) is \cite[Theorem 5.6]{Pey_Sor_2010}, part ii) is \cite[Theorem II.1.]{lions}, part iii) is \cite[Remark 11]{bresis-lions} and part iv) is \cite[Theorem 9]{bresis-lions}.\\

A counterexample for strong convergence in case iv) was given in \cite{Guler}, following the ideas in \cite{Bai}. Therefore, one cannot expect to obtain strong convergence in Theorem \ref{theo2} in general.\\

\subsubsection{Penalized algorithms: diagonal proximal algorithm }
In general, we obtain the diagonal proximal algorithm from \cite{algo1}:
$$\left\{\begin{array}{ll}
x_1 & \in H,\\
x_{n+1} &= (I+\lambda_n A+\lambda_n\beta_n\partial\Psi_2)^{-1}(x_n)\ \ {\rm for\ all}\ n\geq 1,
\end{array}\right.\leqno{(DPA)}$$

Hypothesis ${\rm (H_0)}$ becomes
$$\left\{\begin{array}{rl}
i)&{\mathbf T}=A+N_C\ \hbox{is maximal monotone and}\ S={\mathbf T}^{-1}(0)\neq \emptyset ;\\
ii)& \hbox{For each $z\in N_C(\H)$},\ \sum\limits_{n=1}^\infty
\lambda_n\beta_n\left[\Psi_2^*(\dfrac{z}{\beta_n})-\sigma_C(\dfrac{z}{\beta_n})\right]<\infty;\\
iii)& \sum\limits_{n=1}^\infty\lambda_n=+\infty.
\end{array}\right.\leqno{\rm (H'_0)}$$

Let $(x_n)$ be a sequence generated by (DPA) and define $(z_n)$ and $(\widehat z_n)$ and in \eqref{zn}.

\begin{cor}
Let ${\rm (H'_0)}$ hold. As $n\to\infty$, we have the following:
\begin{itemize}
	\item[i)] The sequence $(\widehat{z}_n)$ converges weakly to a point in $S$;
	\item[ii)] If $(\lambda_n)\in\ell^2$, then the sequence $(z_n)$ converges weakly to a point in $S$; and
	\item[iii)] If $A$ is strongly monotone, then $(x_n)$ converges strongly to the unique point in $S$.
\end{itemize}
\end{cor}

Part i) is \cite[Theorem 3.3]{algo1} and part iii) is \cite[Theorem 3.4]{algo1}. For the weak convergence, we have

\begin{cor}
Let ${\rm (H'_0)}$ hold with $A=\partial\Phi_2$. Assume any of the following conditions holds:
\begin{itemize}
\item[(i)] $\liminf\limits_{n\to\infty}\lambda_n\beta_n>0$ and either $\Phi_2$ or $\Psi_2$ is boundedly inf-compact; or
\item[(ii)] $\liminf\limits_{n\to\infty}\lambda_n\beta_n>0,$ $(\lambda_n)$ is bounded and $\beta_{n+1}-\beta_n\leq K\lambda_{n+1}\beta_{n+1}$
for some $K>0$.
\end{itemize}
Then, the sequence $(x_n)$ converges weakly to a point in $S$ as $n\to\infty$, with $$\lim\limits_{n\to\infty}\Phi_2(x_n)=\min_{C}(\Phi_2).$$
Moreover, convergence is strong in case $(i)$ if $\Psi_2$ is  boundedly inf-compact.
\end{cor}

The hypotheses in part ii) are very close to, but slightly different from those corresponding to cases ii) and iii) of \cite[Theorem 3.6]{algo1}.

\subsection{The forward algorithm}

Taking the operator $A$ and the potential $\Psi_2$ to be zero, the backward step disappears and the  (SFBP) algorithm turns into a purely forward algorithm.

\subsubsection*{The gradient method}

If additionnally $\Psi_1=0$, we obtain the {\em gradient method}, which dates back to~\cite{Cau}:
$$\left\{\begin{array}{rcl}
x_1 & \in & H,\\
x_{n+1} & = & x_n-\lambda_n\nabla\Phi(x_n)\ \ {\rm for\ all}\ n\geq 1.
\end{array}\right.\leqno{(GRAD)}$$

Let $(x_n)$ be a sequence generated by (GRAD).

\begin{cor}
Let $\Phi$ be a convex function with Lipschitz-continuous gradient. Assume $S\neq\emptyset$ and $(\lambda_n)\in\ltwo\setminus\luno$. As $n\to\infty$, the sequence $(x_n)$ converges weakly to a point in $S$, with $\lim\limits_{n\to\infty}\Phi(x_n)=\min(\Phi)$. 
\end{cor}

This is not the most general convergence result for the gradient method. The hypothesis $(\lambda_n)\in\ltwo$ may be replaced by $\limsup_{n\to\infty}\lambda_n<2/L_\Phi$ (see \cite[Teorema 9.6]{Pey_book}). Intermediate results are proved in \cite[Paragraph 1.2.13]{bertsekas}, assuming the step sizes tend to zero; and in \cite[Theorem 3]{Burachik_1995}, under a very precise condition on the step sizes: $\delta_1\leq \lambda_n\leq \frac{2}{L_\Phi}(1-\delta_2)$ with $\delta_1,\delta_1>0$ such that $\frac{L_\Phi}{2} \delta_1+\delta_2<1$. The last two are proved in $\H=\R^n$, but the proof can be easily adapted to the Hilbert-space framework. 

\subsubsection{Penalized algorithm: a diagonal gradient scheme} If $A\equiv 0$ we obtain the {\em diagonal gradient scheme} studied in \cite{algo4}, namely:
$$\ \left\{\begin{array}{rcl}
x_1 & \in & \H,\\
x_{n+1} & = & x_n-\lambda_n\nabla\Phi(x_n)-\lambda_n\beta_n\nabla\Psi_1(x_n),\  n\geq 1.
\end{array}\right.\leqno(DGS)$$

Hypothesis ${\rm (H_0)}$ becomes
$$\left\{\begin{array}{rl}
i)& S={\mathbf T}^{-1}(0)\neq \emptyset ;\\
ii)& \nabla\Phi\ \hbox{is}\ L_{\Phi}\hbox{-Lipschitz-continuous\ and}\ \nabla\Psi_1 \ \hbox{is}\ L_{\Psi_1}\hbox{-Lipschitz-continuous};\\
iii)&\hbox{For each $z\in N_C(\H)$},\ \ \  \sum\limits_{n=1}^\infty
\lambda_n\beta_n\left[\Psi_1^*(\dfrac{z}{\beta_n})-\sigma_C(\dfrac{z}{\beta_n})\right]<\infty;\\
iv)& \sum\limits_{n=1}^\infty\lambda_n=+\infty,\ \ \sum\limits_{n=1}^\infty\lambda^2_n<+\infty,\ \ \hbox{and}\ \
\limsup\limits_{n\to\infty}( L_{\Psi_1}\lambda_n\beta_n)<2.
\end{array}\right.\leqno{\rm (H'''_0)}$$

Then we have:

\begin{cor}
Let ${\rm (H'''_0)}$ hold, and assume that any of the following conditions holds:
\begin{itemize}
\item[(i)] $\liminf\limits_{n\to\infty}\lambda_n\beta_n>0$ and the function $\Phi$ or $\Psi_1$ is boundedly inf-compact; or
\item[(ii)] $\liminf\limits_{n\to\infty}\lambda_n\beta_n>0,$ $(\lambda_n)$ is bounded and $\beta_{n+1}-\beta_n\leq K\lambda_{n+1}\beta_{n+1}$
for some $K>0$.
\end{itemize}
Then, the sequence $(x_n)$ converges weakly to a point in $S$ as $n\to\infty$, with $$\lim\limits_{n\to\infty}\Phi(x_n)=\min_{C}(\Phi).$$
Moreover, convergence is strong in case $(i)$ if $\Psi_1$ is boundedly inf-compact.
\end{cor}

Part ii) was given in \cite[Theorem 2.1]{algo4} with slightly different hypotheses.

\subsection{The forward-backward splitting}

\subsubsection{with no penalization}
If $A\not\equiv 0$ and $\Phi\not\equiv 0$ we obtain the forward-backward splitting:
$$\left\{\begin{array}{rcl}
x_1 & \in & H,\\
x_{n+1} & = & (I+\lambda_n A)^{-1}(x_n-\lambda_n\nabla\Phi(x_n))\ \ {\rm for\ all}\ n\geq 1.
\end{array}\right.\leqno{(FB)}$$

This method combines the previous two, inheriting their virtues and drawbacks. A particularly interesting case is the {\em projected gradient method} described as follows: Let $C$ be a nonempty, closed and convex subset of $\H$ and set $A=N_C$. Then (FB) becomes
$$\left\{\begin{array}{rcl}
x_1 & \in & H,\\
x_{n+1} & = & \hbox{Proj}_C(x_n-\lambda_n\nabla\Phi(x_n))\ \ {\rm for\ all}\ n\geq 1.
\end{array}\right.\leqno{(PG)}$$
This is useful for minimization problems of the form
$$\min\{\Phi(x):x\in C\},$$
when the projection onto the set $C$ is easily performed. \\


Let $(x_n)$ be a sequence generated by (FB).


\begin{cor}
Let $\Phi$ be a convex function with Lipschitz-continuous gradient. Assume $S\neq\emptyset$ and $(\lambda_n)\in\ltwo\setminus\luno$. As $n\to\infty$, we have the following:
\begin{itemize}
	\item[i)] The sequence $(\widehat{z}_n)$ converges weakly to a point in $S$;
	\item[ii)] If $(\lambda_n)\in\ell^2$, then the sequence $(z_n)$ converges weakly to a point in $S$;
	\item[iii)] If $A$ is strongly monotone, then $(x_n)$ converges strongly to the unique point in $S$; and
	\item[iv)] If $A=\partial \Phi_2$, then $(x_n)$ converges weakly to a point in $S$, with $\lim\limits_{n\to\infty}\Phi_2(x_n)=\min(\Phi_2)$.
\end{itemize} 
\end{cor}


The results in \cite[Theorem 1]{Bruck_1977} and \cite[Theorem 2]{Pas_1979} are closely related to part ii). Although they consider a maximal monotone operator $B$ instead of $\nabla\Phi$ (which is a more general framework), their results rely on a $\ltwo$-summability condition $-$ that is difficult to check in practice $-$ concerning a sequence $(w_n)$ satisfying $w_n\in\lambda_nB(x_n)$. Analogously, \cite[Corollary 1]{Pas_1979} is close to part iv).\\

If the step sizes are bounded from below by a positive constant, then the sequence $(x_n)$ converges weakly, even when $A$ is a maximal monotone operator and $\nabla\Phi$ is replaced by a {\em cocoercive} function $B$ (see \cite[Corollary 6.5]{Combettes_2004}). A function $B:\H\to\H$ is cocoercive if $\langle Bx-By,x-y\rangle\ge\beta\|Bx-By\|^2$ for all $x,y\in\H$.

\subsubsection*{Smooth penalization} 
 If $\Phi\equiv 0$, we obtain the forward-backward-penalty scheme studied in \cite{algo3} and \cite{Noun-Pey}:
$$\ \left\{\begin{array}{rcl}
x_1 & \in & \H,\\
x_{n+1} & = & (I+\lambda_n A)^{-1}(x_n-\lambda_n\beta_n\nabla\Psi_1(x_n)),\  n\geq 1.
\end{array}\right.\leqno(FBP)$$

Hypothesis ${\rm (H_0)}$ becomes
$$\left\{\begin{array}{rl}
i)&{\mathbf T}=A+N_C\ \hbox{is maximal monotone and}\ S={\mathbf T}^{-1}(0)\neq \emptyset ;\\
ii)& \nabla\Psi_1 \ \hbox{is}\ L_{\Psi_1}\hbox{-Lipschitz-continuous};\\
iii)&\hbox{For each $z\in N_C(\H)$},\ \ \  \sum\limits_{n=1}^\infty
\lambda_n\beta_n\left[\Psi_1^*(\dfrac{z}{\beta_n})-\sigma_C(\dfrac{z}{\beta_n})\right]<\infty;\\
iv)& \sum\limits_{n=1}^\infty\lambda_n=+\infty,\ \ \hbox{and}\ \
\limsup\limits_{n\to\infty}( L_{\Psi_1}\lambda_n\beta_n)<2.\\
\end{array}\right.\leqno{\rm (H''_0)}$$

\begin{cor}
Assume that ${\rm (H''_0)}$ holds. As $n\to\infty$, we have the following:
\begin{itemize}
	\item[i)] The sequence $(\widehat{z}_n)$ converges weakly to a point in $S$;
	\item[ii)] If $(\lambda_n)\in\ell^2$, then the sequence $(z_n)$ converges weakly to a point in $S$;
	\item[iii)] If $A$ is strongly monotone, then $(x_n)$ converges strongly to the unique point in $S$.
\end{itemize}
\end{cor}

Parts ii) and iii) yield \cite[Theorem 12]{algo3}. Observe that, in part i), convergence is proved without the $\ltwo$-summability assumption, unlike in \cite{algo3}. For the weak convergence we have the following:
\begin{cor}
Let ${\rm (H''_0)}$ hold with $A=\partial \Phi_2$. Assume that any of the following conditions holds:
\begin{itemize}
\item[(i)] $\liminf\limits_{n\to\infty}\lambda_n\beta_n>0$ and the function $\Phi_2$ or $\Psi_1$ is boundedly inf-compact; or
\item[(ii)] $\liminf\limits_{n\to\infty}\lambda_n\beta_n>0,$ $(\lambda_n)$ is bounded and $\beta_{n+1}-\beta_n\leq K\lambda_{n+1}\beta_{n+1}$
for some $K>0$.
\end{itemize}
Then, the sequence $(x_n)$ converges weakly to a point in $S$ as $n\to\infty$, with $$\lim\limits_{n\to\infty}\Phi_2(x_n)=\min_{C}(\Phi_2).$$
Moreover, convergence is strong in case $(i)$ if $\Psi_1$ is boundedly inf-compact.\\
\end{cor}

Part i) yields \cite[Theorem 16]{algo3}, Part ii) yields \cite[Theorem 1]{Noun-Pey}. Both cited results additionally assume the $\ltwo$-summability of the step sizes.

\begin{rem}
The main results of this paper appeared in the PhD thesis of Nahla Noun \cite{Nahla_these}. Simultaneously, Bo\c{t} and Csetnek \cite{Bot_Csetnek} extended the forward-backward results of \cite{algo3} in order to solve the variational inequality
$$0 \in  Ax + Dx + N_C(x),$$
where $A$ is a maximal monotone operator, $D$ a cocoercive operator, and $C$ is the set of zeroes of another maximal monotone operator. Their framework is related but different from ours and cannot be immediately compared.
\end{rem}


\section{An illustration: Sparse-optimal control of a linear system of ODE's} \label{secexamples}Given $y_0\in\R^n$, $A:[0,T]\to\R^{n\times n}$, $B:[0,T]\to\R^{n\times m}$, and $c:[0,T]\to\R^n$, consider the control system
$$\left\{\begin{array}{rcl}
\dot y(t) & = & A(t)y(t)+B(t)u(t)+c(t),\qquad t\in(0,T)\\
y(0) & = & y_0.
\end{array}\right.\leqno{(\cali{CS})}$$
We assume that the functions $A$, $B$ and $c$ are bounded and sufficiently regular so that, for each $u\in L^\infty(0,T;\R^m)$, the system $(\cali{CS})$ has a unique solution $y_u:[0,T]\to\R^n$, which is an absolutely continuous function such that $y_u(0)=y_0$ and the differential equation holds almost everywhere.\\

We are interested in the {\em optimal control problem}
$$\min\,\left\{\,\frac{1}{2}\|y_u-\bar y\|^2_{L^2(0,T;\R^n)}+\|u\|^2_{L^2(0,T;\R^m)}+\|u\|_{L^1(0,T;\R^m)}\,:\,u\in\cali U\,\right\},\leqno{(\cali{OCP})}$$
where $\bar y$ is a {\em reference} trajectory, and the set of admissible controls is
$$\cali U=\{\,u:[0,T]\to\R^m\,:\,\hbox{$u$ is measurable and }|u_i(t)|\le 1\hbox{ a.e. for each }i=1,\dots,m\,\}.$$
The term $\|u\|^2_{L^2(0,T;\R^m)}$ can be interpreted as a measure of the energy invested in controlling the system, and the minimization of the term $\|u\|_{L^1(0,T;\R^m)}$ is known to induce sparsity of the solution.\\

Let $R:[0,T]\to\R^{n\times n}$ be the resolvent of the matrix equation $\dot X=AX$ with initial condition $X(0)=I$. The pair $(u,y)$ satisfies $(\cali CS)$ if, and only if,
$$y(t)=R(t)y_0+R(t)\int_0^tR(s)^{-1}\left[B(s)u(s)+c(s)\right]\,dt.$$
This, in turn, is equivalent to
$$\M(u,y)+z_0=0,$$
where we have written
$$\M(u,y)(t)=-y(t)+R(t)\int_0^tR(s)^{-1}B(s)u(s)\,dt,\  \hbox{ and }\ z_0(t)=R(t)y_0+R(t)\int_0^tR(s)^{-1}c(s)\,dt.$$
Set $H=L^2(0,T;\R^m)\times L^2(0,T;\R^n)$. Since $\M$ is a bounded linear operator from $H$ to $L^2(0,T;\R^n)$, the function $\Psi_1:H\to\R$, defined by
$$\Psi_1(u,y)=\frac{1}{2}\|\M(u,y)+z_0\|^2_{L^2(0,T;\R^n)},$$
is convex and continuously differentiable. On the other hand, since $\cali U$ is nonempty, closed in $L^2(0,T;\R^m)$ and convex, the function $\Psi_2:H\to\Rcupinf$, defined by
$$\Psi_2(u,y)=\delta_{\cali U}(u)$$
(the indicator function of the set $\cali U$), is proper, lower-semicontinuous and convex. Moreover, the pair $(u,y)$ satisfies $(\cali{CS})$ with $u\in\cali U$ if, and only if, $(u,y)\in\argmin(\Psi_1+\Psi_2)$. With this notation, the optimal control problem $(\cali{OCP})$ is equivalent to the {\em constrained optimization problem}
$$\min\,\{\,\Phi_1(u,y)+\Phi_2(u,y)\,:\,(u,y)\in\argmin(\Psi_1+\Psi_2)\,\},\leqno{(\cali{COP})}$$
where $\Phi_1:H\to\R$ is the convex and continuously differentiable function defined by
$$\Phi_1(u,y)=\frac{1}{2}\|y-\bar y\|^2_{L^2(0,T;\R^n)}+\|u\|^2_{L^2(0,T;\R^m)},$$
and $\Phi_2:H\to\Rcupinf$ is the proper, lower-semicontinuous and convex function given by
$$\Phi_2(u,y)=\|u\|_{L^1(0,T;\R^m)}.$$


\section{Proofs}\label{secproposition}

Recall that, by assumption ${\rm (H_0)}$, the monotone operator ${\mathbf T}=A+\nabla \Phi+N_C$ is maximal and $S={\mathbf T}^{-1}(0)\neq \emptyset$. Its domain is $\dom({\mathbf T})=C \cap\dom(A)$. The functions  $\Phi$ and $\Psi_1$ are differentiable and their gradients $\nabla\Phi$ and $\nabla\Psi_1$ are Lipschitz-continuous with constants $L_{\Phi}$ and $L_{\Psi_1}$, respectively. \if{We may assume, without loss of generality, that $\min(\Psi_i)=0$ ($i=1,2$), since otherwise we can use $\widetilde{\Psi}_i=\Psi_i-\min(\Psi_i)$ instead.}\fi Since $\min \Psi_1=\min\Psi_2=0$, the function $\Psi_1+\Psi_2$ vanishes on $C=\argmin(\Psi_1)\cap\argmin(\Psi_2)$.\\

The proofs of Theorems \ref{theo2'}$-$\ref{theo2} ultimately rely on a well-known tool from \cite{Opi} (see the proper statement in \cite{BaiTT}) and \cite{Pas_1979} which gives weak convergence without a priori knowledge of the limit. 

\begin{lem}\label{opial} 
Given a sequence $(x_n)$ in $\H$, a sequence $(\lambda_n)\not\in \ell^1$ of positive numbers, set
$$z_n=\dfrac{1}{\tau_n}\sum\limits_{k=1}^{n}\lambda_kx_k \qquad
\mbox{ and }  \qquad\widehat{z}_n=\dfrac{1}{\tau_n}\sum\limits_{k=1}^{n}\lambda_kx_{k+1},  \qquad \mbox{ with } \qquad \tau_n=\sum\limits_{k=1}^{n}\lambda_k.$$
Let $S$ be a subset of $\H$ and assume that
\begin{itemize}
\item[(i)] for every $x\in S,$  $\lim\limits_{n\to\infty}\|x_n-x\|$ exists;
\item[(ii)]  every weak cluster point of $(x_n),$ respectively $(z_n)$, resp.  $(\widehat{z}_n)$, lies in $S.$
 \end{itemize}
Then $(x_n)$, respectively $(z_n)$, resp. $(\widehat{z}_n)$, converges weakly to a point in $S$ as $n\to\infty$.
\end{lem}

The core of the convergence analysis is the following estimation:

\begin{lem}\label{lem3}
There exist $a,b,c,d,e> 0$ such that, for every $u\in\dom(\mathbf T)$, $z\in Au$, $v\in N_C(u)$, and $w=z+\nabla\Phi(u)+v$, the following inequality holds for  $n$ large enough
\begin{multline}\label{E:lem3}
[1-aL_{\Phi}\lambda_n^2]\|x_{n+1}-u\|^2-\|x_n-u\|^2+b \|x_{n+1}-x_n\|^2+ c\lambda_n\|\nabla\Phi(x_{n+1})-\nabla\Phi(u)\|^2\\
+\dfrac{d}{2}\lambda_n\beta_n(\Psi_1+\Psi_2)(x_{n+1})+e\lambda_n\beta_n\|\nabla\Psi_1(x_n)\|^2\\ \leq\dfrac{d}{2}\lambda_n\beta_n\left[(\Psi_1+\Psi_2)^*(\frac{4v}{d\beta_n})-\sigma_C(\frac{4v}{d\beta_n})\right]
+2\lambda_n\langle w,u- x_{n+1}\rangle.
\end{multline}
\end{lem}

The proof uses standard convex analysis tools along with very careful estimations. Since it is highly technical, it will be given below in Subsection \ref{SS:Proof_Lemma}. A straightforward consequence of Lemma \ref{lem3} is the following proposition, which contains the basic properties of the algorithm, including the first assumption of Lemma \ref{opial}:

\begin{prop} \label{prop1} Assume ${\rm (H_0)}$, $(\lambda_n L_\Phi)\in \ell^2$, and let $(x_n)$ be a sequence generated by the (SFBP) Algorithm.
Then the following holds:
\begin{itemize}
\item [i)] For every $u\in S$, $\lim\limits_{n\to\infty}\|x_n-u\|$ exists.
\item [ii)] The series $\sum\limits_{n\geq1}\|x_{n+1}-x_n\|^2$, $\sum\limits_{n\geq1}\lambda_n\beta_n\Psi_1(x_n)$, $\sum\limits_{n\geq1}\lambda_n\beta_n\Psi_2(x_n),$ $\sum\limits_{n\geq1}\lambda_n\|\nabla\Phi(x_{n+1})-\nabla \Phi(u)\|^2$ and $\sum\limits_{n\geq1}\lambda_n\beta_n\|\nabla\Psi_1(x_{n})\|^2$ converge.
\item [iii)] If moreover $\liminf\limits_{n\to\infty}\lambda_n\beta_n>0,$ then $\sum\limits_{n\geq1}\Psi_1(x_n)$ and $\sum\limits_{n\geq1}\Psi_2(x_n)$ converge, $\lim\limits_{n\to\infty}\Psi_1(x_n)=\lim\limits_{n\to\infty}\Psi_2(x_n)=0,$ and every weak cluster point of the sequence $(x_n)$ lies in $C$.\\
\end{itemize}
\end{prop}

Indeed, if $u\in S$, then $0\in A u+\nabla\Phi(u)+N_C(u)$. Write $0=z+\nabla\Phi(u)+v$ with $z\in Au$ and $v\in  N_C(u)$. For $n$ large enough, Lemma~\ref{lem3} gives
\begin{multline*}
[1-aL_{\Phi}\lambda_n^2]\|x_{n+1}-u\|^2-\|x_n-u\|^2+b \|x_{n+1}-x_n\|^2+ c\lambda_n\|\nabla\Phi(x_{n+1})-\nabla\Phi(u)\|^2\\
+\dfrac{d}{2}\lambda_n\beta_n(\Psi_1+\Psi_2)(x_{n+1})+e\lambda_n\beta_n\|\nabla\Psi_1(x_n)\|^2
 \leq\dfrac{d}{2}\lambda_n\beta_n\left[(\Psi_1+\Psi_2)^*(\frac{4v}{d\beta_n})-\sigma_C(\frac{4v}{d\beta_n})\right]
\end{multline*}
Since the right-hand side is summable, all the parts of Proposition \ref{prop1} ensue from the following elementary fact concerning the convergence of real sequences:

\begin{lem}\label{lemfondamental}
Let $(a_n)$, $(\delta_n)$ and $(\varepsilon_n)$ be nonnegative and let $(\xi_n)$ be bounded from below. Assume
$$ 
(1-a_n)\xi_{n+1}-\xi_n+\delta_n\le \varepsilon_n
$$ 
for all $n$ large enough. If $(a_n)$ and $(\varepsilon_n)$ belong to $\luno$, then $(\xi_n)$ is convergent and $(\delta_n)$ belongs to $\luno$.
\end{lem}

\if{
\dem We may assume, without loss of generality, that $(\xi_n)$ is nonnegative (otherwise, we subtract any lower bound). Inequality \eqref{E:1} can be rewritten as
$$\xi_{n+1}-\xi_n+\delta_n\le \varepsilon_n+a_n\xi_{n+1}.$$
If $(\xi_n)$ is bounded, then the right-hand side is in $\luno$ and everything is a direct consequence of the classical convergence result \cite[Lemma 2]{algo3}. Since $(\xi_n)$ is nonnegative, the proof is complete if we prove that it is bounded from above.
We may also assume that $1-a_n>0$ because $\limn a_n=0$. Since $\delta_n\ge 0$, \eqref{E:1} gives 
$$\xi_{n+1}\le\frac{1}{1-a_n}\left[\xi_n+\varepsilon_n\right].$$
Observing that $(1-a_n)^{-1}\le\prod_{k=1}^n(1-a_k)^{-1}$, a simple induction argument shows that
$$\xi_{n+1}\le\left[\prod_{k=1}^n(1-a_k)^{-1}\right]\left[\xi_1+\sum_{k=1}^n\varepsilon_k\right].$$
The right-hand side is bounded because $(a_n)$ and $(\varepsilon_n)$ belong to $\luno$. \phantom{ksjdh}\bx
}\fi


\subsection{Proof of Lemma \ref{lem3}}\label{SS:Proof_Lemma} Take $u\in\dom(\mathbf T)$, $z\in Au$, $v\in N_C(u)$, and $w=z+\nabla\Phi(u)+v$. Rewrite algorithm (SFBP) as
\begin{equation}\label{algo6}
\frac{x_n-x_{n+1}}{\lambda_n}-\nabla\Phi(x_n)-v_{n+1}-\beta_nw_{n+1}-\beta_n\nabla\Psi_1(x_n)=0.
\end{equation}
with $v_{n+1}\in Ax_{n+1} $ and $w_{n+1}\in\partial\Psi_2(x_{n+1})$.

\begin{claim}\label{claim1}
The following inequality holds for every $n$:
\begin{multline} \label{equ1}
\|x_{n+1}-u\|^2-\|x_n-u\|^2+\|x_{n+1}-x_n\|^2
+2\lambda_n\langle\nabla\Phi(x_n)-\nabla\Phi(u),x_{n+1}-u\rangle
 \\+2\lambda_n\beta_n\langle \nabla\Psi_1(x_n),x_{n+1}-u\rangle
 \leq
2\lambda_n\langle \nabla\Phi(u)+z,u-x_{n+1}\rangle-2\lambda_n\beta_n\Psi_2(x_{n+1}).
\end{multline}
\end{claim}

\dem The monotonicity of $A$ at points $u$ and $x_{n+1}$ gives
\begin{equation}\label{Amono}
2\lambda_n\langle v_{n+1},u- x_{n+1}\rangle \leq 2\lambda_n\langle z,u- x_{n+1}\rangle.
\end{equation}
The subdifferential inequality for function $\Psi_2$ at $x_{n+1}$ with $w_{n+1}\in \partial\Psi_2(x_{n+1})$ gives
\begin{equation}\label{psi2}
2\lambda_n\beta_n\langle w_{n+1},u- x_{n+1}\rangle\leq -2\lambda_n\beta_n\Psi_2(x_{n+1}).
\end{equation}
To conlude, it suffices to sum \eqref{Amono} and \eqref{psi2}, use \eqref{algo6}, along with the fact that
$$
2\langle x_{n+1}-x_n, x_{n+1}-u\rangle =\|x_{n+1}-u\|^2-\|x_n-u\|^2+\|x_{n+1}-x_n\|^2,
$$ 
and rearrange the terms conveniently. \bx

The next step is to estimate some of the terms in \eqref{equ1} (observe that the last two terms on the left-hand side vanish if $\nabla\Phi$ is constant or $\Psi_1\equiv 0$, which correspond to the cases $L_{\Phi}=0$ and $L_{\Psi_1}=0$, respectively). At different points we shall use the Baillon-Haddad Theorem \cite{baillon-haddad}: 

\begin{lem}[Baillon-Haddad Theorem]\label{lembaillon}
Let $f:\H\rightarrow \R$ be a convex differentiable function and  $L_f>0$. Then $\nabla f$ is Lipschitz continuous with constant $L_f$ if and only if
 $\nabla f$ is $\frac{1}{L_f}$-cocoercive.
\end{lem}

We also use the following Descent Lemma (see, for example \cite{bertsekas}):

\begin{lem}[Descent Lemma]\label{lemdescent}
Let $f:\H\rightarrow \R$ be continuously differentiable such that $\nabla f$ is Lipschitz continuous with constant $L_f$. Then, for every  $x$ and $y$  in $\H$,
$$f(x+y)\leq f(x) + \langle\nabla f(x),y\rangle + \dfrac{L_f}{2}\|y\|^2.$$
\end{lem}

We have the following:

\begin{claim}\label{claim2}
Assume $\nabla\Phi$ is not constant. For every $\eta>0$ we have:
$$2\lambda_n\langle\nabla\Phi(x_n)-\nabla\Phi(u),x_{n+1}-u\rangle\geq \dfrac{-\lambda_n^2}{\eta}\|x_{n+1}-u\|^2-\eta L_\Phi^2\|x_{n+1}-x_n\|^2 +\dfrac{2\lambda_n}{L_\Phi}\|\nabla\Phi(x_{n+1})-\nabla\Phi(u)\|^2.$$
\end{claim}

\dem Write
$$ 
\langle\nabla\Phi(x_n)-\nabla\Phi(u),x_{n+1}-u\rangle= \langle\nabla\Phi(x_n)-\nabla\Phi(x_{n+1}),x_{n+1}-u\rangle+\langle\nabla\Phi(x_{n+1})-\nabla\Phi(u),x_{n+1}-u\rangle.
$$ 
We easily see that
$$ 
2\lambda_n\langle\nabla\Phi(x_n)-\nabla\Phi(x_{n+1}),x_{n+1}-u\rangle\geq\dfrac{-1}{\eta} \lambda_n^2\|x_{n+1}-u\|^2-\eta L_\Phi^2\|x_{n+1}-x_n\|^2.
$$ 
On the other hand, Lemma \ref{lembaillon} implies
$$ 
2\lambda_n\langle\nabla\Phi(x_{n+1})-\nabla\Phi(u),x_{n+1}-u\rangle\geq \dfrac{2\lambda_n}{L_\Phi}\|\nabla\Phi(x_{n+1})-\nabla\Phi(u)\|^2.
$$ 
The result follows immediately.\bx


\begin{claim}\label{claim3}
Assume ${\Psi_1}\not\equiv 0$. For all $\eta,\theta>0$ and $n\in \N$ we have
\begin{multline*} 
2\lambda_n\beta_n\langle \nabla\Psi_1(x_n),x_{n+1}-u\rangle\geqslant \left[\dfrac{2}{(1+\eta)L_{\Psi_1}}-\dfrac{(1+\theta)}{1+\eta}\lambda_n\beta_n\right]\lambda_n\beta_n\|\nabla\Psi_1(x_n)\|^2\\
+\dfrac{2\eta}{1+\eta}\lambda_n\beta_n\Psi_1(x_{n+1})-
\left[\dfrac{1}{(1+\theta)(1+\eta)}+\dfrac{\eta L_{\Psi_1}}{1+\eta}\lambda_n\beta_n\right]\|x_{n+1}-x_n\|^2.
\end{multline*}
\end{claim}

\dem Write
\begin{equation}\label{equ9'}
2\lambda_n\beta_n\langle\nabla\Psi_1(x_n),x_{n+1}-u\rangle=
2\lambda_n\beta_n\langle\nabla\Psi_1(x_n),x_{n+1}-x_n\rangle
+2\lambda_n\beta_n\langle\nabla\Psi_1(x_n),x_{n}-u\rangle.
\end{equation}
A convex combination of the bounds given by Lemma \ref{lembaillon} and the subdifferential inequality gives
\begin{equation}\label{e3}
\langle\nabla\Psi_1(x_n),x_{n}-u\rangle\geq \dfrac{1}{(1+\eta)L_{\Psi_1}}\|\nabla\Psi_1(x_n)\|^2
+\dfrac{\eta}{1+\eta}\Psi_1(x_n)
\end{equation}
for any $\eta>0$. Now take $\theta>0$ and use the identity
\begin{multline*}
\dfrac{1}{1+\theta}\|x_{n+1}-x_n+(1+\theta)\lambda_n\beta_n\nabla\Psi_1(x_n)\|^2=\\
\dfrac{1}{1+\theta}\|x_{n+1}-x_n\|^2+(1+\theta)\lambda_n^2\beta_n^2\|\nabla\Psi_1(x_n)\|^2
+2\lambda_n\beta_n\langle\nabla\Psi_1(x_n),x_{n+1}-x_n\rangle,
 \end{multline*}
to obtain
\begin{equation*} 
2\lambda_n\beta_n\langle\nabla\Psi_1(x_n),x_{n+1}-x_n\rangle\geq -\dfrac{1}{1+\theta}\|x_{n+1}-x_n\|^2-(1+\theta)\lambda_n^2\beta_n^2\|\nabla\Psi_1(x_n)\|^2.
\end{equation*}
On the other hand, Lemma \ref{lemdescent} at $x_n$ and $x_{n+1}$ gives
$$2\lambda_n\beta_n\langle\nabla\Psi_1(x_n),x_{n+1}-x_n\rangle\geq
2\lambda_n\beta_n[\Psi_1(x_{n+1})-\Psi_1(x_n)]
-L_{\Psi_1}\lambda_n\beta_n\|x_{n+1}-x_n\|^2.$$
A convex combination of the last two inequalities produces
\begin{multline}\label{e5}
2\lambda_n\beta_n\langle\nabla\Psi_1(x_n),x_{n+1}-x_n\rangle\geq
-\dfrac{1+\theta}{1+\eta}\lambda_n^2\beta_n^2\|\nabla\Psi_1(x_n)\|^2\\
+\dfrac{2\eta}{1+\eta}\lambda_n\beta_n[\Psi_1(x_{n+1})-\Psi_1(x_n)]
-\left[\dfrac{1}{(1+\eta)(1+\theta)}+\dfrac{\eta L_{\Psi_1}}{1+\eta}\lambda_n\beta_n\right]\|x_{n+1}-x_n\|^2.
\end{multline}
Finally, use \eqref{e3} and \eqref{e5} in \eqref{equ9'} to conclude.\bx


\begin{claim}\label{chap2-lem2}
There exist $a,b,c,d,e>0$ such that for all sufficiently large $n\in \N$ we have
\begin{multline}\label{eqlem2}
\left(1-aL_{\Phi}\lambda_n^2\right)\|x_{n+1}-u\|^2-\|x_n-u\|^2+b \|x_{n+1}-x_n\|^2
+ c\lambda_n\|\nabla\Phi(x_{n+1})-\nabla\Phi(u)\|^2\\+d\lambda_n\beta_n(\Psi_1+\Psi_2)(x_{n+1})+e\lambda_n\beta_n\|\nabla\Psi_1(x_n)\|^2
 \leq
2\lambda_n\langle \nabla\Phi(u)+z,u-x_{n+1}\rangle.
\end{multline}
\end{claim}

\dem We focus on the case where $\nabla \Phi$ is not constant and $\Psi_1\not\equiv 0$. The other cases are simpler and left to the reader. Claims \ref{claim1}, \ref{claim2} and \ref{claim3}, and the fact that
$$-2\lambda_n\beta_n\Psi_2(x_{n+1})\leq \dfrac{-2\eta}{1+\eta}\lambda_n\beta_n\Psi_2(x_{n+1})$$
for every $\eta>0$, together imply
\begin{multline*}
\left[1-\frac{\lambda_n^2}{\eta}\right]\|x_{n+1}-u\|^2-\|x_{n}-u\|^2+ \left [1-\frac{1}{(1+\theta)(1+\eta)}
-\frac{\eta L_{\Psi_1}}{1+\eta}\lambda_n\beta_n-\eta L_{\Phi}^2\right]\|x_{n+1}-x_n\|^2\\
+\dfrac{2}{L_{\Phi}}\lambda_n\|\nabla\Phi(x_{n+1})-\nabla\Phi(u)\|^2
+\left[\dfrac{2}{(1+\eta)L_{\Psi_1}}-\dfrac{(1+\theta)}{1+\eta}\lambda_n\beta_n\right]\lambda_n\beta_n\|\nabla\Psi_1(x_n)\|^2\\
+\frac{2\eta}{1+\eta}\lambda_n\beta_n(\Psi_1+\Psi_2)(x_{n+1})\leq
2\lambda_n\langle \nabla\Phi(u)+z,u-x_{n+1}\rangle
\end{multline*}
for every $\eta,\theta>0$. Set $\Gamma=\limsup\limits_{n\to\infty}L_{\Psi_1}\lambda_n\beta_n<2$ and take $\theta_0>0$ small enough such that
$$2-(1+\theta_0)\lambda_n\beta_n L_{\Psi_1}\geq 2-(1+\theta_0)\Gamma >0$$
for all sufficiently large $n$. Since
$$
\lim\limits_{\eta\to 0^+}\left[1-\dfrac{1}{(1+\theta)(1+\eta)}-\dfrac{\eta\Gamma L_{\Psi_1}}{1+\eta}-\eta L_\Phi^2\right]=1-\dfrac{1}{(1+\theta)}>0,
$$
we can take $\eta_0>0$ so that \eqref{eqlem2} holds with $a=\dfrac{1}{\eta_0L_{\Phi}}$, $b=1-\dfrac{1}{(1+\theta_0)(1+\eta_0)}-\dfrac{\eta_0\Gamma L_{\Psi_1}}{1+\eta_0}-\eta_0 L_\Phi^2$, $c=\dfrac{2}{L_\Phi},\quad d=\dfrac{2\eta_0}{1+\eta_0}$, and $e=\dfrac{2(1+\eta_0)(1+\theta_0)L_{\Psi_1}\Gamma}{(1+\eta_0)L_{\Psi_1}}$, which are all positive.\bx

\noindent{\bf Proof of Lemma \ref{lem3}, completed.} Observe that
\begin{align*}
2\lambda_n\langle \nabla\Phi(u)+z,\ & u-x_{n+1}\rangle - \dfrac{d}{2}\lambda_n\beta_n(\Psi_1+\Psi_2)(x_{n+1}) \\
&=2\lambda_n\langle v, x_{n+1}-u\rangle - \dfrac{d}{2}\lambda_n\beta_n(\Psi_1+\Psi_2)(x_{n+1}) -2\lambda_n\langle w, x_{n+1}-u\rangle\\
&=\dfrac{d}{2} \lambda_n\beta_n\left[\langle \dfrac{4v}{d\beta_n}, x_{n+1}\rangle-(\Psi_1+\Psi_2)(x_{n+1})-\langle\dfrac{4v}{d\beta_n}, u\rangle\right] -2\lambda_n\langle w, x_{n+1}-u\rangle\\
&\leq\dfrac{d}{2}\lambda_n\beta_n\left[(\Psi_1+\Psi_2)^*\left(\frac{4v}{d\beta_n}\right)-\sigma_C\left(\frac{4v}{d\beta_n}\right)\right]-2\lambda_n\langle w, x_{n+1}-u\rangle
\end{align*}
because $\dfrac{4v}{d\beta_n}\in N_C(u)$ implies $\sigma_C\left(\dfrac{4v}{d\beta_n}\right)=\langle \dfrac{4v}{d\beta_n}, u\rangle$. Whence
\begin{multline}\label{equ14}
2\lambda_n\langle p-z, x_{n+1}-u\rangle\leq \dfrac{d}{2}\lambda_n\beta_n(\Psi_1+\Psi_2)(x_{n+1})\\
+\dfrac{d}{2}\lambda_n\beta_n\left[(\Psi_1+\Psi_2)^*\left(\frac{4v}{d\beta_n}\right)-\sigma_C\left(\frac{4v}{d\beta_n}\right)\right]
+2\lambda_n\langle w, u-x_{n+1}\rangle.
\end{multline}
We obtain \eqref{E:lem3} by using \eqref{equ14} in \eqref{eqlem2} and rearranging the terms containing $(\Psi_1+\Psi_2)(x_{n+1})$.\bx


\subsection{Weak ergodic convergence: proof of Theorem \ref{theo2'}}\label{secergodic}\label{S:ergodic}

In view of Lemma \ref{opial} and part i) of Proposition \ref{prop1}, it suffices to prove that every weak cluster point of the sequence $(z_n)$, respectively $(\widehat{z}_n)$, lies in $S$. By maximal monotonicity of ${\mathbf T}$, a point $\overline x$ belongs to $S$ if and only if $\langle w,u-\overline x\rangle\geq 0$ for all $u \in C\cap \hbox{dom}(A)$ and all $w \in {\mathbf T}(u)$.\\

We begin with $(\widehat{z}_n)$. Take any $u\in\ C\ \cap$ dom$(A)$ and $w\in {\mathbf T}(u)$. By Lemma \ref{lem3}, we have
\begin{multline}\label{equ15}
\|x_{n+1}-u\|^2-\|x_n-u\|^2
 \leq
\dfrac{d}{2}\lambda_n\beta_n\left[(\Psi_1+\Psi_2)^*(\frac{4v}{d\beta_n})-\sigma_C(\frac{4v}{d\beta_n})\right]\\
+aL_{\Phi}\lambda_n^2\|x_{n+1}-u\|^2
+2\lambda_n\langle w,u-x_{n+1}\rangle
\end{multline}
for $n$ large enough. Since $\|x_{n+1}-u\|$ converges as $n\to \infty$, it is bounded. Let $a\|x_{n+1}-u\|^2\leq M$ for some $M >0$ and every $n$.
Take $$
\varepsilon_n=\dfrac{d}{2}\lambda_n\beta_n\left[(\Psi_1+\Psi_2)^*(\frac{4v}{d\beta_n})-\sigma_C(\frac{4v}{d\beta_n})\right]
+ML_{\Phi}\lambda_n^2.
$$
Assumption ${\rm (H_0)}\ iii)$ and $(\lambda_nL_{\Phi})\in \ell^2$ yield $\sum\limits_{n\geq1}\varepsilon_n < +\infty$.
Summing up for $k=1,...,n,$ we have
$$
\|x_{n+1}-u\|^2-\|x_1-u\|^2
 \leq
2\langle w,\sum\limits_{k=1}^n\lambda_ku\rangle
-2\langle w,\sum\limits_{k=1}^n\lambda_kx_{k+1}\rangle
+\sum\limits_{k=1}^n\varepsilon_k.
$$
Removing the nonnegative term $\|x_{n+1}-u\|^2$ and dividing by $2\tau_n=2\sum\limits_{k=1}^n\lambda_k,$
we get
\begin{equation}\label{equ17}
\dfrac{-\|x_1-u\|^2}{2\tau_n} \leq \langle w,u-\widehat{z}_n\rangle + \dfrac{1}{2\tau_n}\sum\limits_{k=1}^n\varepsilon_k
\end{equation}
Passing to the lower limit in (\ref{equ17}) and using $\tau_n\rightarrow \infty$ as $n\to\infty$ (because $\lambda_n \not\in \luno$)
we deduce that
$$
\liminf\limits_{n\to\infty}\langle w,u-\widehat{z}_n\rangle\geq 0.
$$
If some subsequence $(\widehat{z}_{n_k})$ converges weakly to $x_\infty,$ then $\langle w,u-x_\infty\rangle\geq 0$. Thus $x_\infty\in S$.
\bigskip

For the sequence $(z_{n})$, we decompose the term $\langle w,u-x_{n+1}\rangle$ in \eqref{equ15} and write
\begin{multline*}
\|x_{n+1}-u\|^2-\|x_n-u\|^2
 \leq
\dfrac{d}{2}\lambda_n\beta_n\left[(\Psi_1+\Psi_2)^*(\frac{4v}{d\beta_n})-\sigma_C(\frac{4v}{d\beta_n})\right]\\
+aL_{\Phi}\lambda_n^2\|x_{n+1}-u\|^2
+2\lambda_n\langle w,x_n-x_{n+1}\rangle+2\lambda_n\langle w,u-x_n\rangle.
\end{multline*}
Using
$2\lambda_n\langle w,x_n-x_{n+1}\rangle \leq \lambda_n^2\|w\|^2 + \|x_{n+1}-x_n\|^2$
in the last inequality and proceeding as above we obtain
$$
\|x_{n+1}-u\|^2-\|x_1-u\|^2
 \leq
2\langle w,\sum\limits_{k=1}^n\lambda_ku\rangle
-2\langle w,\sum\limits_{k=1}^n\lambda_kx_k\rangle
+\sum\limits_{k=1}^n\zeta_k,
$$
where
$$
\zeta_n=\dfrac{d}{2}\lambda_n\beta_n\left[(\Psi_1+\Psi_2)^*(\frac{4v}{d\beta_n})-\sigma_C(\frac{4v}{d\beta_n})\right]
+(M L_{\Phi}+\|w\|^2)\lambda_n^2+\|x_{n+1}-x_n\|^2.
$$
Assumption ${\rm (H_0)}\ iii),$ Proposition \ref{prop1} $ii)$ and the additional assumption $(\lambda_n)\in \ell^2$ give
$\sum\limits_{n\geq1}\zeta_n < +\infty$ and we conclude as before.\bx

\subsection{Strong convergence: proof of Theorem \ref{theo1'}}\label{secstrong}\label{S:strong}

In this Section we treat the particular case where additionally the maximal monotone operator $A$ is strongly monotone, or the potential $\Phi$ is strongly convex. Recall that $A$ is strongly monotone with parameter $\alpha>0$ if
$$\langle x^*-y^*,x-y\rangle \geq \alpha \|x-y\|^2$$
whenever $x^*\in Ax$ and $y^*\in Ay$.  The potential $\Phi$ is strongly convex if $\nabla \Phi$ is  strongly monotone. Since $A$, $\nabla \Phi$ and $N_C$ are monotone, the operator ${\mathbf T}$ is also strongly monotone whenever $A$  is strongly monotone or $\Phi$ is strongly convex.
Then, the set ${\mathbf T}^{-1}0$ reduces to the singleton $\{u\}$, for some $u\in \H$. By the definition of $S$, there exist $z\in Au$ and $v\in N_C(u)$ such that $z+\nabla\Phi(u)+v=0$.\\

The proof of Theorem \ref{theo1'} is a direct consequence of the following reinforced version of Lemma \ref{lem3}:

\begin{lem}
There exist $a,b,c,d,e >0$ such that, for $n$ large enough, we have
\begin{multline*}
\alpha\lambda_n\|x_{n+1}-u\|^2+\|x_{n+1}-u\|^2-\|x_n-u\|^2+b\|x_{n+1}-x_n\|^2
+ c\lambda_n\|\nabla\Phi(x_{n+1})-\nabla\Phi(u)\|^2\\
+\dfrac{d}{2}\lambda_n\beta_n(\Psi_1+\Psi_2)(x_{n+1})+e\lambda_n\beta_n\|\nabla\Psi_1(x_n)\|^2\\
 \leq
\dfrac{d}{2}\lambda_n\beta_n\left[(\Psi_1+\Psi_2)^*(\frac{4v}{d\beta_n})-\sigma_C(\frac{4v}{d\beta_n})\right]
+aL_{\Phi}\lambda_n^2\|x_{n+1}-u\|^2.
\end{multline*}
\end{lem}


\subsection{Weak convergence: proof of Theorem \ref{theo2}}\label{secweak}\label{S:weak}
This section achieves the proof of  Theorem \ref{theo2}, that is the  weak convergence of the sequence $(x_n)$ generated by the  (SFBP) algorithm,
 in the special case where $A=\partial\Phi_2$ is the subdifferential of a proper lower-semicontinuous
convex function $\Phi_2:\H\rightarrow \R\cup \{+\infty\}.$ Writing  $\Phi_1$ instead of  $\Phi$,  for the sake of symmetry, the (SFBP) algorithm takes the form
\begin{equation}\label{sous_dif}
\left\{\begin{array}{ll}
x_1 & \in \H,\\
x_{n+1} &= (I+\lambda_n \partial\Phi_2+\lambda_n\beta_n\partial\Psi_2)^{-1}(x_n-\lambda_n\nabla\Phi_1(x_n)-\lambda_n\beta_n\nabla\Psi_1(x_n))\ \forall\ n\geq 1.
\end{array}\right.\hfill
\end{equation}
Since $\partial\Phi_2+\nabla\Phi_1+N_C$ is maximal monotone, the solution set $S$ is equal to
$$
S=(\partial\Phi_2+\nabla\Phi_1+N_C)^{-1}(0)=
\argmin\{\Phi_1(u)+\Phi_2(u):u\in \argmin\Psi_1\cap\argmin \Psi_2\}.
$$

\noindent

We prove the weak convergence of the sequence $(x_n)$ generated by   algorithm (\ref{sous_dif})  to some point in $S$ using Opial-Passty's Lemma~\ref{opial}. The first assumption  is satisfied from Proposition \ref{prop1} i). The second assumption,  that every weak cluster point of $(x_n)$ belongs to $S$, will be verified in three different cases $(i)$, respectively $(ii)$ and $(iii)$, in Subsection \ref{subsecinf}, respectively Subsection \ref{subsecbeta}.
\paragraph{Finally,} in Subsection \ref{subsecprooftheo} we finish the proof of Theorem \ref{theo2}  with the minimizing property of the sequence.

\subsubsection{Weak convergence in case $(i)$:  bounded inf-compactness}\label{subsecinf}

Denote by dist(.,$S$) the distance function to the closed convex set $S$ and set $d(x)=\frac{1}{2}$dist$(x,S)^2.$
The function $d$ is convex and differentiable, $\nabla d(x)=x-P_S(x)$ where $P_S$ denotes the projection onto $S.$

The proof goes along the same lines as that of {\cite[Theorem 16]{algo3}}. 
In the next lemma, we prove that  $\lim\limits_{n\to\infty}d(x_n)=0.$ By the weak lower semicontinuity of the convex function $d$, it implies that every weak cluster point of $(x_n)$ lies in $S.$ Thus $(x_n)$ satisfies the second assumption of  Opial-Passty's Lemma, and we deduce the weak convergence to some point in $S$.
Besides, since  $\liminf\limits_{n\to\infty}\lambda_n\beta_n>0$ by assumption, Proposition \ref{prop1} iii) yields  $\lim\limits_{n\to\infty}(\Psi_1+\Psi_2)(x_n)=0$. If we additionnaly assume that $(\Psi_1+\Psi_2)$ is  boundedly inf-compact, the bounded sequence $(x_n)$ is also relatively compact. Hence its  weak convergence implies its strong convergence.
 This achieves the proof of the first part of Theorem \ref{theo2} in Case (i).

\begin{lem} Under the assumptions of Theorem \ref{theo2} (i), let $(x_n)$ be a sequence generated by the (SFBP) algorithm. Then
$$\lim\limits_{n\to\infty}d(x_n)=0.$$
\end{lem}

\dem We reformulate (\ref{sous_dif}) as
\begin{equation}\label{0}
x_n-x_{n+1}=\lambda_n\nabla\Phi_1(x_n)+\lambda_n\beta_n\nabla\Psi_1(x_n)+\lambda_nv_{n+1}+\lambda_n\beta_nw_{n+1},
\end{equation}
where $v_{n+1}\in\partial\Phi_2(x_{n+1}) $ and $w_{n+1}\in\partial\Psi_2(x_{n+1}).$
The convexity of $d$ with (\ref{0}) yields
\begin{eqnarray}\label{1}
d(x_n)&\geq & d(x_{n+1}) +\langle x_{n+1}-P_S(x_{n+1}),x_{n}-x_{n+1}\rangle\nonumber\\
&=& d(x_{n+1})+\lambda_n\langle \nabla\Phi_1(x_n)+v_{n+1},x_{n+1}-P_S(x_{n+1})\rangle\\
&& +\lambda_n\beta_n\langle \nabla\Psi_1(x_n),x_{n+1}-P_S(x_{n+1})\rangle
+\lambda_n\beta_n\langle w_{n+1},x_{n+1}-P_S(x_{n+1})\rangle.\nonumber
\end{eqnarray}
We treat each term on the right-hand side of (\ref{1}).
Let
$$
\alpha=\hbox{min}\{(\Phi_1+\Phi_2)(z);\ z\in C\}.
$$
Firstly for $\lambda_n\langle \nabla\Phi_1(x_n)+v_{n+1},x_{n+1}-P_S(x_{n+1})\rangle.$ Since $\Phi_1$ is convex we have
\begin{eqnarray}\label{2}
\Phi_1(P_Sx_{n+1})
&\geq& \Phi_1(x_{n})+ \langle \nabla\Phi_1(x_n),P_Sx_{n+1}-x_n\rangle\\
&=&\Phi_1(x_{n})+ \langle \nabla\Phi_1(x_n),P_Sx_{n+1}-x_{n+1}\rangle +\langle \nabla\Phi_1(x_n),x_{n+1}-x_n\rangle.\nonumber
\end{eqnarray}
From Descent Lemma (Lemma \ref{lemdescent}) we have
$$
\Phi_1(x_{n+1})\leq  \Phi_1(x_{n})+ \langle \nabla\Phi_1(x_n),x_{n+1}-x_n\rangle+\frac{L_{\Phi_1}}{2}\|x_{n+1}-x_n\|^2.
$$
Using this in (\ref{2}) it follows that
\begin{equation}\label{3}
\Phi_1(P_Sx_{n+1})\geq \Phi_1(x_{n+1})- \frac{L_{\Phi_1}}{2}\|x_{n+1}-x_n\|^2 + \langle \nabla\Phi_1(x_n),P_Sx_{n+1}-x_{n+1}\rangle.
\end{equation}
On the other hand, the subdifferential inequality for $\Phi_2$ writes
\begin{equation}\label{4}
\Phi_2(P_Sx_{n+1})\geq \Phi_2(x_{n+1})+ \langle v_{n+1},P_Sx_{n+1}-x_{n+1}\rangle.
\end{equation}
Noting that $\Phi_1(P_Sx_{n+1})+\Phi_2(P_Sx_{n+1})=\alpha$ and adding (\ref{3}) and (\ref{4}) we get
\begin{equation}\label{5}
\langle \nabla\Phi_1(x_n)+v_{n+1},P_Sx_{n+1}-x_{n+1}\rangle \leq \alpha- (\Phi_1+\Phi_2)(x_{n+1})+\frac{L_{\Phi_1}}{2}\|x_{n+1}-x_n\|^2.
\end{equation}
Secondly for $\lambda_n\beta_n\langle \nabla\Psi_1(x_n),x_{n+1}-P_S(x_{n+1})\rangle.$ If $\Psi_1=0$ then it is equal to zero. If $\Psi_1\neq0,$ since $\nabla\Psi_1(P_Sx_{n+1})=0,$ the cocoercivity of $\nabla\Psi_1$ implies
$$
\langle \nabla\Psi_1(x_n),P_Sx_{n+1}-x_{n}\rangle \leq -\frac{1}{L_{\Psi_1}}\|\nabla\Psi_1(x_n)\|^2.
$$
Adding the last inequality and
$$\langle \nabla\Psi_1(x_n),x_{n}-x_{n+1}\rangle\leq \frac{1}{L_{\Psi_1}}\|\nabla\Psi_1(x_n)\|^2+\frac{L_{\Psi_1}}{4}\|x_{n+1}-x_n\|^2$$
 we deduce
\begin{equation}\label{6}
\langle \nabla\Psi_1(x_n),P_Sx_{n+1}-x_{n+1}\rangle
\leq \frac{L_{\Psi_1}}{4}\|x_{n+1}-x_n\|^2.
\end{equation}
Thirdly for $\lambda_n\beta_n\langle w_{n+1},x_{n+1}-P_S(x_{n+1})\rangle.$ Since $w_{n+1}\in \partial\Psi_2(w_{n+1})$ and $0\in \partial\Psi_2(P_Sx_{n+1}),$
the monotonicity of $\partial\Psi_2$ implies
\begin{equation}\label{7}
\langle w_{n+1},x_{n+1}-P_Sx_{n+1}\rangle \geq 0.
\end{equation}
Combining (\ref{5}), (\ref{6}) and (\ref{7}) in (\ref{1}),  and since $\limsup\limits_{n\to\infty}( L_{\Psi_1}\lambda_n\beta_n)<2$,
 we deduce that
\begin{eqnarray*}
d(x_{n+1}) -d(x_{n})+\lambda_n[(\Phi_1+\Phi_2)(x_{n+1})-\alpha]&\leq& \lambda_n\frac{L_{\Phi_1}}{2}\|x_{n+1}-x_n\|^2+\lambda_n\beta_n\frac{L_{\Psi_1}}{4}\|x_{n+1}-x_n\|^2\\
\label{8} &\leq&\frac{1}{2}\|x_{n+1}-x_n\|^2,
\end{eqnarray*}
for $n$ large enough.
The remainder of the proof is an adaptation of the proof of {\cite[Theorem 16]{algo3}}: considering that $x_{n+1}$ may not lie in $C$, and we may not have $(\Phi_1+\Phi_2)(x_{n+1})-\alpha\ge 0$ for every $n\in\N$, it is achieved by studying separately the two cases:\\ 
\noi Case I: There exists $n_0\in\N$ such that $(\Phi_1+\Phi_2)(x_n)\ge \alpha$ for all $n\ge n_0$.\\
\noi Case II: For each $n\in\N$ there exists $n'>n$ such that $(\Phi_1+\Phi_2)(x_{n'})<\alpha$.

\if{

\noindent
The right-hand side of inequality (\ref{8}) is summable by Part $ii)$ in Proposition \ref{prop1}. But $x_{n+1}$ may not lie in $C$, and we may not have $(\Phi_1+\Phi_2)(x_{n+1})-\alpha\ge 0$ for every $n\in\N$.
 We shall analyze the two possible situations separately, namely\\

\noi Case I: There exists $n_0\in\N$ such that $(\Phi_1+\Phi_2)(x_n)\ge \alpha$ for all $n\ge n_0$.\\
\noi Case II: For each $n\in\N$ there exists $n'>n$ such that $(\Phi_1+\Phi_2)(x_{n'})<\alpha$.\\

In our analysis we follow the arguments in \cite{cabot}, \cite{AC} and \cite{algo3} which can be traced back to \cite{Baillon-Cominetti}.\\

\noi\u{Case I:} For every $n\ge n_0$ we have
$$d(x_{n+1})-d(x_n)\le \frac{M}{2}\|x_{n+1}-x_n\|^2.$$
Since $d(x_n)$ is bounded from below and the right-hand side is summable, we conclude that $\limn d(x_n)$ exists. In order to verify that this limit must be $0$ it suffices to find a subsequence that converges to $0$. Let us sum up inequality (\ref{8}) for $n=n_0,\dots,N$ to obtain
$$d(x_{N+1})-d(x_{n_0})+\sum_{n=n_0}^N\lambda_n\left[(\Phi_1+\Phi_2)(x_{n+1})-\alpha\right]
\le\frac{M}{2}\sum_{n=n_0}^N\|x_{n+1}-x_n\|^2.$$
Letting $N\to\infty$ we deduce that $\sum\limits_{n=n_0}^\infty\lambda_n[(\Phi_1+\Phi_2)(x_{n+1})-\alpha]<+\infty$. Since $(\lambda_n)\notin\luno$ we must have $\liminf\limits_{n\to\infty}(\Phi_1+\Phi_2)(x_n)\le \alpha$. Consider a subsequence $(x_{k_n})$ such that
$$\limn(\Phi_1+\Phi_2)(x_{k_n})=\liminf\limits_{n\to\infty}(\Phi_1+\Phi_2)(x_n).$$
Clearly the sequence $(\Phi_1+\Phi_2)(x_{k_n})$ is bounded. By Proposition \ref{prop1} $i)$ , the sequence $(x_n)$ is bounded and so $(\Psi_1+\Psi_2)(x_{k_n})$ is bounded as well.
The  bounded inf-compactness assumption of $(\Phi_1+\Phi_2)$ or $(\Psi_1+\Psi_2)$ ensures the existence of subsequence, that we still denote by  $(x_{k_n})$ without any loss of generality,
that converges strongly to some $\bar x$. By Proposition \ref{prop1} $iii)$ $\bar x$ must belong to $C$. Then $\alpha\leq (\Phi_1+\Phi_2)(\bar x).$
From the weak lower-semicontinuity of $(\Phi_1+\Phi_2)$ we deduce that
$$\alpha\le(\Phi_1+\Phi_2)(\bar x)\le\limn(\Phi_1+\Phi_2)(x_{k_n'})=\liminf_{n\to\infty}(\Phi_1+\Phi_2)(x_n)\le \alpha.$$
This shows that $\bar x\in\S.$  View that  $(x_{k_n'})\rightarrow \bar x$ and $d$ is continuous we deduce $\limn d(x_{k_n'})=0$.\\

\noi\u{Case II:} For all sufficiently large $n$ the number
$$\tau_n=\max\{k\le n:(\Phi_1+\Phi_2)(x_k)<\alpha\}$$
is well defined. Observe that we have $\limn\tau_n=+\infty$.
Take $N\in\N$ (large enough for $\tau_N$ to exist).
If $\tau_N<N$ then $(\Phi_1+\Phi_2)(x_{n+1})\ge \alpha$ for $n=\tau_N,\dots,N-1$.
Inequality (\ref{8}) then gives
$$d(x_N)-d(x_{\tau_N})\le\frac{M}{2}\sum_{n=\tau_N}^{N-1}\|x_{n+1}-x_n\|^2.$$
If $N=\tau_N$ then $d(x_{\tau_N})=d(x_N)$. In either case we have
$$d(x_N)-d(x_{\tau_N})\le\frac{M}{2}\sum_{n=\tau_N}^{\infty}\|x_{n+1}-x_n\|^2$$
and so letting $N\to\infty$ we deduce that
$$0\le\limsup_{n\to\infty}d(x_n)\le\limsup_{n\to\infty}d(x_{\tau_n}).$$
It suffices to prove that $\limsup\limits_{n\to\infty}d(x_{\tau_n})=0$.
Since $(\Phi_1+\Phi_2)(x_{\tau_n})<\alpha$ for all $n$ one has $\limsup\limits_{n\to\infty}(\Phi_1+\Phi_2)(x_{\tau_n})\le \alpha$.
In particular the sequence $(\Phi_1+\Phi_2)(x_{\tau_n})$ is bounded.
As before, the same is true for $(\Psi_1+\Psi_2)(x_{\tau_n})$ and using the  bounded inf-compactness assumption one concludes,
as in case I, that a subsequence converges strongly to a point in $\S$, which guarantees that $\limn d(x_n)=0$.
}\fi

\bx

\subsubsection{Weak convergence in cases $(ii)$ and $(iii)$: bounded increase of the sequence $(\beta_n)$ and the unconstrained case}\label{subsecbeta}
As before, it suffices to  prove that every weak cluster point of the sequence $(x_n)$ generated by the (SFBP) algorithm lies in $S$ to deduce 
its  weak convergence to a point in $S$. We decompose the proof in several lemmas.
Let us introduce the penalized functions $\Omega_n$ and $\widetilde{\Omega}_n$ defined on $\H$ by
$$\Omega_n=\Phi_1+\beta_n\Psi_1 \hspace*{1cm} \mbox{and} \hspace*{1cm} \widetilde{\Omega}_n=\Phi_2+\beta_n\Psi_2.$$
Being the sum of two smooth functions whose gradient is lipschitz-continuous, $\Omega_n$ is in his turn a smooth function
 whose gradient is lipschitz-continuous with constant $L_n=L_{\Phi_1}+\beta_nL_{\Psi_1}$. If $\Psi_1=0,$  $\Omega_n$ reduces to $\Phi_1$ and $L_n=L_{\Phi_1}.$

\begin{lem}\label{chap2-lem5}
Assume that ${\rm (H_0)}$ hold with $A=\partial \Phi_2$, $(\lambda_nL_{\Phi})\in \ell^2$ and let $(x_n)$ be a sequence generated by the $(SFBP)$ algorithm. Then the following holds:
\begin{itemize}
\item[i)] For every $n\geq 1$ the penalized functions verify:
\begin{multline*}
\left[\Omega_{n+1}(x_{n+1})+\widetilde{\Omega}_{n+1}(x_{n+1})\right]-\left[\Omega_n(x_n)+\widetilde{\Omega}_n(x_n)\right]
\\+\left[\frac{1}{\lambda_n}-\frac{L_n}{2}\right]\|x_{n+1}-x_n\|^2
\le
(\beta_{n+1}-\beta_n)(\Psi_1+\Psi_2)(x_{n+1}).
\end{multline*}
\item[ii)] If~~~~ $\beta_{n+1}-\beta_n\leq K\lambda_{n+1}\beta_{n+1}$, or if  $\Psi_1=\Psi_2=0$,  then the sequence $\left(\Omega_n(x_n)+\widetilde{\Omega}_n(x_n)\right)$ converges as $n\to \infty$.
\end{itemize}
\end{lem}
{\bf Proof of Part i).}
Apply Descent Lemma~\ref{lemdescent} to obtain
$$
\Omega_n(x_{n+1})-\Omega_n(x_n)\leq \langle \nabla\Omega_n(x_n),x_{n+1}-x_n\rangle+\frac{L_n}{2}\|x_{n+1}-x_n\|^2.
$$
Using $\Omega_n(x_{n+1})=\Omega_{n+1}(x_{n+1})-(\beta_{n+1}-\beta_n)\Psi_1(x_{n+1})$, it follows that
\begin{equation}
\Omega_{n+1}(x_{n+1})-\Omega_n(x_n)\leq \langle \nabla\Omega_n(x_n),x_{n+1}-x_n\rangle+\frac{L_n}{2}\|x_{n+1}-x_n\|^2+(\beta_{n+1}-\beta_n)\Psi_1(x_{n+1}).\label{eq1}
\end{equation}
On the other hand, remark that
\begin{equation}
\label{eqomegatild}
\widetilde{\Omega}_{n+1}(x_{n+1}) -\widetilde{\Omega}_n(x_n)
= \Phi_2(x_{n+1})-\Phi_2(x_n) + \beta_n\left[\Psi_2(x_{n+1})-\Psi_2(x_n)\right] + (\beta_{n+1}-\beta_n)\Psi_2(x_{n+1}).
\end{equation}
Recall the following formulation of the (SFBP) algorithm
\begin{equation*}
x_n-x_{n+1}=\lambda_n\nabla\Phi_1(x_n)+\lambda_n\beta_n\nabla\Psi_1(x_n)+\lambda_nv_{n+1}+\lambda_n\beta_nw_{n+1},
\end{equation*}
with $v_{n+1}\in\partial\Phi_2(x_{n+1}) $ and $w_{n+1}\in\partial\Psi_2(x_{n+1}).$
The subdifferential inequality at $x_{n+1}$ of $\Phi_2$ and $\Psi_2,$ respectively gives
$$
\Phi_2(x_{n+1})-\Phi_2(x_n) \leq \langle v_{n+1}, x_{n+1}-x_n\rangle
\hspace{0.5cm}{\rm and}\hspace{0.5cm}
\Psi_2(x_{n+1})-\Psi_2(x_n) \leq \langle w_{n+1}, x_{n+1}-x_n\rangle.
$$
Replacing this in (\ref{eqomegatild}) we obtain
$$
\widetilde{\Omega}_{n+1}(x_{n+1})-\widetilde{\Omega}_n(x_n)\leq \langle v_{n+1}+\beta_n w_{n+1}, x_{n+1}-x_n\rangle + (\beta_{n+1}-\beta_n)\Psi_2(x_{n+1}).
$$
Adding (\ref{eq1}) and the last inequality we deduce that
\begin{multline}\label{eq7}
\left[\Omega_{n+1}(x_{n+1}) + \widetilde{\Omega}_{n+1}(x_{n+1})\right]-\left[\Omega_n(x_n)+\widetilde{\Omega}_n(x_n)\right]\leq 
 \langle \nabla\Omega_n(x_n)+v_{n+1}+\beta_n w_{n+1},x_{n+1}-x_n\rangle\\
 +\frac{L_n}{2}\|x_{n+1}-x_n\|^2 +(\beta_{n+1}-\beta_n)(\Psi_1+\Psi_2)(x_{n+1}).
\end{multline}
Therefore, just substitute the equality $\nabla\Omega_n(x_n)+v_{n+1}+\beta_n w_{n+1}=-\frac{x_{n+1}-x_n}{\lambda_n}$ to conclude $i)$. \\

{\bf Proof of Part ii).} Since $L_{\Psi_1}\lambda_n\beta_n <2$ for $n$ large enough, from $i)$ we have
\begin{multline}\label{130624}
\left[\Omega_{n+1}(x_{n+1})+\widetilde{\Omega}_{n+1}(x_{n+1})\right]-\left[\Omega_n(x_n)+\widetilde{\Omega}_n(x_n)\right]\\
\le (\beta_{n+1}-\beta_n)(\Psi_1+\Psi_2)(x_{n+1})+\dfrac{L_{\Phi_1}}{2}\|x_{n+1}-x_n\|^2.
\end{multline}
Take an element $u\in S$ and $z\in \partial \Phi_2(u)$. Write the subdifferential inequality at $u$ for $\Phi_1$ and $\Phi_2$ to obtain
 \begin{equation}\begin{split}
\Omega_n(x_n)+\widetilde{\Omega}_n(x_n)
&\geq \Phi_1(x_n)+\Phi_2(x_n)\\
&\geq\Phi_1(u)+\Phi_2(u)+\langle \nabla\Phi_1(u)+z, x_n-u\rangle.\label{eq}
\end{split}
\end{equation}
Since  $(x_n)$ is bounded by Proposition \ref{prop1} i), the sequence $\left(\Omega_n(x_n)+\widetilde{\Omega}_n(x_n)\right)$ is bounded from below.
Now, if $\beta_{n+1}-\beta_n\leq K\lambda_{n+1}\beta_{n+1}$, or  if $\Psi_1=\Psi_2=0$, the right-hand side of (\ref{130624}) is summable from  Proposition \ref{prop1} ii). This implies the convergence of the sequence $\left(\Omega_n(x_n)+\widetilde{\Omega}_n(x_n)\right)$.\bx

\begin{lem}\label{prop2}  Assume ${\rm (H_0)}$, $(\lambda_nL_{\Phi})\in \ell^2$, and additionnally $\lambda_n$ is bounded and $\beta_{n+1}-\beta_n\leq K\lambda_{n+1}\beta_{n+1}$,  respectively $\Psi_1=\Psi_2=0$. Let $(x_n)$ be a sequence generated by the (SFBP) algorithm. Then, for every $u\in S$ we have
 \begin{equation}\label{equ22}
\sum\limits_{n\geq1}\lambda_n\left[\Omega_{n+1}(x_{n+1})+\widetilde{\Omega}_{n+1}(x_{n+1})-\Phi_1(u)-\Phi_2(u)\right]<+\infty \ \ \ (\rm{possibly -\infty}).
\end{equation}
\end{lem}
\dem We write
\begin{multline}\label{equ24}
2\lambda_n\left[\Omega_{n+1}(x_{n+1})+\widetilde{\Omega}_{n+1}(x_{n+1})-\Phi_1(u)-\Phi_2(u)\right]=\\
2\lambda_n\left[\Phi_1(x_{n+1})-\Phi_1(u)\right]+2\lambda_n\left[\Phi_2(x_{n+1})-\Phi_2(u)\right]
+2\lambda_n\beta_{n+1} (\Psi_1+\Psi_2)(x_{n+1}).
\end{multline}

The (SFBP) algorithm writes
\begin{equation*}
v_{n+1}=\dfrac{x_n-x_{n+1}}{\lambda_n}-\nabla\Phi_1(x_n)-\beta_n\nabla\Psi_1(x_n)-\beta_n w_{n+1}
\end{equation*}
with $v_{n+1}\in\partial\Phi_2(x_{n+1}) $ and $w_{n+1}\in\partial\Psi_2(x_{n+1}).$
The subdifferential inequality of $\Phi_1$ and $\Phi_2,$ respectively gives
$$
\Phi_1(x_{n+1})-\Phi_1(u)\leq \langle \nabla\Phi_1(x_{n+1}),x_{n+1}-u\rangle
\hspace{0.5cm}{\rm and}\hspace{0.5cm}
\Phi_2(x_{n+1})-\Phi_2(u)\leq \langle v_{n+1},x_{n+1}-u\rangle.
$$
Thus
\begin{multline*}
2\lambda_n\left[\Omega_{n+1}(x_{n+1})+\widetilde{\Omega}_{n+1}(x_{n+1})-\Phi_1(u)-\Phi_2(u)\right]\leq \\
2\lambda_n\langle \nabla \Phi_1(x_{n+1}),x_{n+1}-u\rangle +2\lambda_n\langle v_{n+1}, x_{n+1}-u\rangle
+2\lambda_n\beta_{n+1} (\Psi_1+\Psi_2)(x_{n+1}).
\end{multline*}
Write $2\langle x_{n+1}-x_n, x_{n+1}-u\rangle=\|x_{n+1}-x_n\|^2+\|x_{n+1}-u\|^2-\|x_n-u\|^2$ and deduce
\begin{multline}\label{equ27}
2\lambda_n\left[\Omega_{n+1}(x_{n+1})+\widetilde{\Omega}_{n+1}(x_{n+1})-\Phi_1(u)-\Phi_2(u)\right]
\leq -\|x_{n+1}-x_n\|^2-\|x_{n+1}-u\|^2+\|x_n-u\|^2\\
+2\lambda_n\langle\nabla\Phi_1( x_{n+1})-\nabla\Phi_1(x_n), x_{n+1}-u\rangle
+2\lambda_n\beta_n\langle\nabla\Psi_1( x_n),u- x_{n+1}\rangle\\
+2\lambda_n\beta_n\langle w_{n+1},u- x_{n+1}\rangle+2\lambda_n\beta_{n+1}\Psi_2(x_{n+1})
+2\lambda_n\beta_{n+1}\Psi_1(x_{n+1}).
\end{multline}

We now treat each term on the right-hand side of (\ref{equ27}).
\paragraph{For the term} $2\lambda_n\langle\nabla\Phi_1( x_{n+1})-\nabla\Phi_1(x_n), x_{n+1}-u\rangle$ on the right-hand side of (\ref{equ27}), use the Cauchy-Schwartz inequality to write
\begin{equation*}
\begin{split}
2\lambda_n\langle\nabla\Phi_1( x_{n+1})-\nabla\Phi_1(x_n), x_{n+1}-u\rangle
&\leq2\lambda_n\|\nabla\Phi_1( x_{n+1})-\nabla\Phi_1(x_n)\| \|x_{n+1}-u\|\\
&\leq2 \lambda_n  L_{\Phi_1} \|x_{n+1}-x_n\|\|x_{n+1}-u\|\\
&\leq  L_{\Phi_1} \|x_{n+1}-x_n\|^2+\lambda_n^2L_{\Phi_1}\|x_{n+1}-u\|^2
\end{split}
\end{equation*}
Since $\sum\limits_{n\geq1}\|x_{n+1}-x_n\|^2<+\infty$, $\sum\limits_{n\geq1}\lambda_n^2L_{\Phi_1}<+\infty$
and $(\|x_n-u\|)$ is bounded, we deduce that
\begin{equation}\label{equ28}
\sum\limits_{n\geq1}2\lambda_n\langle\nabla\Phi_1( x_{n+1})-\nabla\Phi_1(x_n), x_{n+1}-u\rangle<+\infty.
\end{equation}

\paragraph{For the term} $2\lambda_n\beta_n\langle\nabla\Psi_1( x_n),u- x_{n+1}\rangle$, 
write
\begin{equation}\label{equ29}
2\lambda_n\beta_n\langle\nabla\Psi_1( x_n),u- x_{n+1}\rangle=2\lambda_n\beta_n\langle\nabla\Psi_1( x_n),u- x_n\rangle
+2\lambda_n\beta_n\langle\nabla\Psi_1( x_n),x_n- x_{n+1}\rangle.
\end{equation}
On one hand, the monotonicity of the gradient and the fact $u\in C$ imply
\begin{equation}\label{equ30}
\langle\nabla\Psi_1( x_n),u- x_n\rangle\leq 0.
\end{equation}
On the other hand we have
\begin{equation*}
2\lambda_n\beta_n\langle\nabla\Psi_1( x_n),x_n- x_{n+1}\rangle\leq
\lambda_n^2\beta_n^2\|\nabla\Psi_1( x_n)\|^2+\|x_{n+1}-x_n\|^2.
\end{equation*}
Therefore Proposition \ref{prop1} ii) and the bound $\lambda_n\beta_n<\frac{2}{L_{\Psi_1}}$,  if $L_{\Psi_1}\ne 0$ and for $n$ large enough, yield
\begin{equation}\label{equ31'}
\sum\limits_{n\geq1}\lambda_n\beta_n \langle\nabla\Psi_1( x_n),x_n- x_{n+1}\rangle<+\infty.
\end{equation}
Combining (\ref{equ30}) and (\ref{equ31'}) in (\ref{equ29}) we conclude
\begin{equation}\label{equ32}
\sum\limits_{n\geq1}\lambda_n\beta_n\langle\nabla\Psi_1( x_n),u- x_{n+1}\rangle < +\infty.
\end{equation}
If  $L_{\Psi_1}= 0$, use the Cauchy-Schwartz inequality and write
\begin{equation*}
2\lambda_n\beta_n\langle\nabla\Psi_1( x_n),u- x_{n+1}\rangle\leq 2\lambda_n\beta_n \|\nabla\Psi_1( x_n)\| \|u- x_{n+1}\|,
\end{equation*}
which is summable, since  $\|u- x_{n+1}\|$ is bounded, $\|\nabla\Psi_1( x_n)\|=a \|\nabla\Psi_1( x_n)\|^2$  with $a>0$, and in view of Proposition \ref{prop1} ii).

\paragraph{For the term} $2\lambda_n\beta_n\langle w_{n+1},u- x_{n+1}\rangle+2\lambda_n\beta_{n+1}\Psi_2(x_{n+1}),$ if $\Psi_2=0$, its value is zero. Otherwise, we assume $\beta_{n+1}-\beta_n\leq K \lambda_{n+1}\beta_{n+1}$ and
use the subdifferential inequality for $\Psi_2$ at $x_{n+1}$ and $u$
$$
\langle w_{n+1},u- x_{n+1}\rangle\leq -\Psi_2(x_{n+1})
$$
to write
 $$2\lambda_n\beta_n\langle w_{n+1},u- x_{n+1}\rangle+2\lambda_n\beta_{n+1}\Psi_2(x_{n+1})
 \leq2K\lambda_n \lambda_{n+1}\beta_{n+1}\Psi_2(x_{n+1}).$$
 Since $\lambda_n$ is bounded by assumption and $\sum\limits_{n\geq1}\lambda_{n+1}\beta_{n+1}\Psi_2(x_{n+1})<\infty$ by Proposition \ref{prop1} $iii)$,
 it follows that
%
 \begin{equation}\label{equ33}
\sum\limits_{n\geq1}[ \lambda_n\beta_n\langle w_{n+1},u- x_{n+1}\rangle+
\lambda_n\beta_{n+1}\Psi_2( x_{n+1})] < +\infty.
\end{equation}
\paragraph{For the remaining term} $2\lambda_n\beta_{n+1}\Psi_1(x_{n+1})$, it is equal to zero if $\Psi_1=0$.  Otherwise, we assume $\beta_{n+1}-\beta_n\leq K \lambda_{n+1}\beta_{n+1}$, and we write
$$
\lambda_n\beta_{n+1}\Psi_1(x_{n+1})=\lambda_n\beta_n\Psi_1(x_{n+1})+ \lambda_n (\beta_{n+1}-\beta_n)\Psi_1(x_{n+1}).
$$
Noting that,  $\lambda_n\beta_n<2/L_{\Psi_1}$ for $n$  large enough, $\lambda_n$ is bounded by assumption  and that $\lambda_n (\beta_{n+1}-\beta_n)\leq K\lambda_n\lambda_{n+1}\beta_{n+1}$,  Proposition \ref{prop1} $iii)$ yields
$\sum\limits_{n\geq1} \lambda_n\beta_n\Psi_1(x_{n+1})<+\infty$ and $\sum\limits_{n\geq1} \lambda_n(\beta_{n+1}-\beta_n)\Psi_1(x_{n+1})<+\infty.$
We then deduce that
\begin{equation}\label{equ34}
\sum\limits_{n\geq1} \lambda_n\beta_{n+1}\Psi_1(x_{n+1})<+\infty.
\end{equation}
Finally we conclude (\ref{equ22}) from (\ref{equ27}) by using (\ref{equ28}), (\ref{equ32}), (\ref{equ33}), (\ref{equ34}) and the fact that
\begin{equation}\label{ee}
\sum\limits_{n\geq1}\|x_n-u\|^2-\|x_{n+1}-u\|^2\leq\|x_1-u\|^2.
\end{equation}
\bx

\begin{lem}\label{lem7}
Assume ${\rm (H_0)}$ and $(\lambda_nL_{\Phi})\in \ell^2$. Assume moreover that one of the following conditions holds:
\begin{itemize}
\item[(i)] $\liminf\limits_{n\to\infty}\lambda_n\beta_n>0,$ $(\lambda_n)$ is bounded and $\beta_{n+1}-\beta_n\leq K\lambda_{n+1}\beta_{n+1}$
for some $K>0$.
\item[(ii)] $\Psi_1=\Psi_2=0$.
\end{itemize}
Let $(x_{n_k})$ be a subsequence of $(x_n)$
that converges weakly to some $x_\infty$ as $n\to\infty$.
Then
$$x_\infty\ \in\ S={\rm Argmin} \{\Phi_1(x)+\Phi_2(x):\ x\in {\rm Argmin}(\Psi_1+\Psi_2)\}. $$
\end{lem}
\dem  Since $\sum\limits_{n\geq1}\lambda_n=+\infty$ by the last statement of Hypotheses ${\rm (H_0)}$,
Lemmas \ref{chap2-lem5} $ii)$ and  \ref{prop2} together imply
\begin{equation}\label{eq1'}
\lim\limits_{n\to\infty}\left[\Omega_n(x_n)+\widetilde{\Omega}_n(x_n)\right] \leq \Phi_1(u)+\Phi_2(u) \ \ \ \forall u \in S.
\end{equation}
Now, in view of (\ref{eq1'}), the weak lower-semicontinuity of $\Phi_1$ and $\Phi_2$ yields
\begin{align}\label{semicont}
\Phi_1(x_\infty) + \Phi_2(x_\infty)
&\leq \liminf\limits_{k\to\infty}\Phi_1(x_{n_k})+\liminf\limits_{k\to\infty}\Phi_2(x_{n_k})\\
&\leq \liminf\limits_{k\to\infty}\Omega_{n_k}(x_{n_k})+\liminf\limits_{k\to\infty}\widetilde{\Omega}_{n_k}(x_{n_k})\nonumber\\
&\leq \liminf\limits_{k\to\infty}\left[\Omega_{n_k}(x_{n_k})+\widetilde{\Omega}_{n_k}(x_{n_k})\right]\nonumber\\
&\qquad = \lim\limits_{n\to\infty}\left[\Omega_n(x_n)+\widetilde{\Omega}_n(x_n)\right]\nonumber\\
&\leq \Phi_1(u)+\Phi_2(u) \ \ \ \ \ \forall u\in S.\nonumber
\end{align}
Under Assumption (i),  $x_\infty$  belongs to $C$ by Proposition \ref{prop1} $iii)$. Under Assumption (ii), $C=\H$. Thus $x_\infty\in S$ and  every weak cluster point of $(x_n)$ lies in $S$.\bx\\



\subsubsection{Minimizing property in cases $(ii)$ and $(iii)$}\label{subsecprooftheo} The first part of  Theorem  \ref{theo2}, the weak convergence of the sequence $(x_n)$,  is  proved in Subsection \ref{subsecinf} in Case $(i)$, and in Subsection \ref{subsecbeta} in Cases $(ii)$ and $(iii)$.
For the second part, recalling that $\nabla\Phi_1(u)=-p$ and $z\in\partial \Phi_2(u)$, the subdifferential inequality of $ \Phi_1+\Phi_2$ at $u\in S$ yields
\begin{equation}\label{liminf}
(\Phi_1 + \Phi_2)(x_n)
\geq (\Phi_1 + \Phi_2)(u)+\langle z-p, x_n-u\rangle.
\end{equation}
Passing to the lower limit in (\ref{liminf}), using  $p-z\in N_C(u)$ and the fact that $(x_n)$ weakly converges to a point in $S$ (First part of Theorem \ref{theo2}), it follows that
$$
\liminf\limits_{n\to\infty}(\Phi_1 + \Phi_2)(x_n)
\geq (\Phi_1 + \Phi_2)(u).
$$
On the other hand, using (\ref{eq1'}) we have
$$
\limsup\limits_{n\to\infty}(\Phi_1 + \Phi_2)(x_n)
\leq \lim\limits_{n\to\infty}\left[\Omega_n(x_n)+\widetilde{\Omega}_n(x_n)\right]
 \leq \Phi_1(u)+\Phi_2(u).
$$
Combining the last two inequalities with the fact that $u\in S$, the second part of Theorem \ref{theo2} directly follows. 
\bx



\begin{thebibliography}{00}

\bibitem{Alv_Pey_2010} F. \'Alvarez and J. Peypouquet, {\em Asymptotic almost-equivalence of Lipschitz evolution systems in Banach spaces}, Nonlinear Anal. 73 (2010), no. 9, 3018--3033.

\bibitem{AttCom} H. Attouch and R. Cominetti, {\em A dynamical approach to convex minimization coupling approximation with the steepest descent method}, J. Differential Equations, 128 (1996), 519--540.

\bibitem{AC} H. Attouch and M.-O. Czarnecki, {\em Asymptotic behavior of coupled dynamical systems with multiscale aspects}, J. Differ. Equations 248 (2010) no. 6, 1315--1344.

\bibitem{algo1} H. Attouch, M.-O. Czarnecki and J. Peypouquet, {\em Prox-penalization and splitting methods for constrained variational problems}, SIAM J. Optim., 21 (2011) no. 1, 149--173

\bibitem{algo3} H. Attouch, M.-O. Czarnecki and J. Peypouquet, {\em Coupling forward-backward with penalty schemes and parallel splitting for constrained
variational inequalities}, SIAM J. Optim. 21 (2011) no. 4, 1251--1274.


\bibitem{Bai} J.-B. Baillon, {\em Un exemple concernant le comportement asymptotique de la solution du probl\`eme $du/dt+\partial\phi(u)=0$}, J. Functional Anal. 28, (1978), 369--376.


\bibitem{BaiTT} J.-B. Baillon, Comportement asymptotique des contractions et semi-groupes de contractions - equations de schroedinger non lineaires et divers, Th\`ese, Paris VI, 1978.


\bibitem{BaiBre} J.-B. Baillon and H. Br\'ezis, {\em Une remarque sur le comportement asymptotique des semi-groupes
non lin\'eaires}, Houston J. Math. 2, (1976), 5--7.


\bibitem{Baillon-Cominetti} J.-B. Baillon and R. Cominetti, {\em A convergence result for nonautonomous subgradient evolution equations and its application to the steepest descent exponential penalty trajectory in linear programming}, J. Funct. Anal. 187 (2001), 263--273.


\bibitem{baillon-haddad} J.-B. Baillon and G. Haddad, {\em Quelques propri\'et\'es des op\'erateurs angle-born\'es et n-cycliquement monotones}, Israel J. Math. 26 (1977), no. 2, 137--150.


\bibitem{bertsekas} D. Bertsekas, Nonlinear programming. Athena Scientific, Belmont MA, 2nd Printing (2003).


\bibitem{Bot_Csetnek} R. I. Bo\c{t}, E. R. Csetnek {\em Forward-Backward and Tseng's type penalty schemes for monotone inclusion problems}, Set-Valued and Variational Analysis 22(2) (2014), 313--331 

\bibitem{brezis} H. Br\'ezis, Op\'erateurs maximaux monotones et semi-groupes de contractions dans les espaces de Hlibert, North Holland Publishing
Company, Amsterdam, 1973.

\bibitem{bresis-lions} H. Br\'ezis and P.-L. Lions, {\em Produits infinis de r\'esolvantes}, Israel J. Math. 29 (1978), 329--345.

\bibitem{Bruck} R.-E. Bruck, {\em Asymptotic convergence of nonlinear contraction semigroups in Hilbert spaces},
 J. Funct. Anal. 18 (1975), 15--26.

\bibitem{Bruck_1977} R.-E. Bruck, {\em On the weak convergence of an ergodic iteration for the solution of variational inequalities for monotone operators in Hilbert space}, J. Math. Anal. Appl. 61 (1977), 159--164.

\bibitem{Burachik_1995} R. Burachik, L.-M. Grana Drummond, A.-N. Iseum and B.-F. Svaiter, {\em Full convergence of the steepest Descent method with inexact line searches}, Optimization 32 (1995), 137--146.

\bibitem{Cab1} A. Cabot, {\em The steepest descent dynamical system with control. Applications to constrained minimization}, ESAIM Control Optim. Calc. Var. 10 (2004), 243--258.

\bibitem{cabot} A. Cabot, {\em Proximal point algorithm controlled by a slowly vanishing term. Applications to hierarchical minimization}, SIAM J.
Optim. 15 (2005), no. 8, 1207--1223.

\bibitem{Cau} A.-L. Cauchy, {M\'ethode g\'en\'erale pour la r\'esolution des syst\`emes d'\'equations simultan\'ees}, C. R. Acad. Sci. Paris 25 (1847), 536--538.

\bibitem{Combettes_2004} P.-L. Combettes, {\em Solving monotone inclusions via compositions of nonexpansive averaged operators}, Optimization 53 (2004), no. 5-6, 475--504.




\bibitem{Eva} L.-C. Evans, Partial differential equations, Second edition, Graduate Studies in Mathematics 19, American Mathematical Society, Providence, RI, 2010.

\bibitem{Guler} O. G\"uler, {\em On the convergence of the proximal point algorithm for convex minimization}, SIAM J. Control Optim. 29 (1991), no. 2, 403--419.

\bibitem{lions} P.-L. Lions, {\em Une methode iterative de resolution d'une inequation variationnelle}, Israel J. Math. 31 (1978), no. 2, 204--208.

\bibitem{Mar1} B. Martinet, {\em R\'egularisation d'in\'equations variationelles par approximations successives}, RAIRO 4 (1970), 154--159.

\bibitem{Mar2} B. Martinet, {\em D\'etermination approch\'ee d'un point fixe d'une application pseudo-contractante}, C.R. Acad. Sci. Paris. 274 (1972), 163--165.

\bibitem{Moreau} J.-J. Moreau, {\em Proximit\'e et dualit\'e dans un espace hilbertien}, Bull. Soc. Math. France 93 (1965), 273-?299.

\bibitem{Nahla_these} N.Noun, {\em Convergence et stabilisation de syst\`emes dynamiques coupl\'es et multi-\'echelles vers des \'equilibres sous contraintes ; application \'a l'optimisation hi\'erarchique}, PhD Thesis, Universit\'e Montpellier 2 and  Universit\'e Libanaise \`a Beyrouth, 2013.


\bibitem{Noun-Pey} N. Noun and J. Peypouquet, {\em Forward-Backward-Penalty scheme for constrained convex minimization without inf-compactness},  J. Optim. Theory Appl. 158 (2013), no. 3, 787--795.

\bibitem{Opi} Z. Opial, {\em Weak Convergence of the sequence of successive approximations for nonexpansive mappings}, Bull. Amer. Math. Soc. 73 (1967), 591--597.

\bibitem{Pas_1979} G. Passty, {\em Ergodic convergence to a zero of the sum of monotone operators in Hilbert space}, J. Math. Anal. Appl. 72 (1979), no. 2, 383--390.

\bibitem{Pey_book} J. Peypouquet, Optimizaci\'on y sistemas din\'amicos, Ediciones IVIC, Caracas, 2013.

\bibitem{algo4} J. Peypouquet, {\em Coupling the gradient method with a general exterior penalization scheme for convex minimization}, J. Optim. Theory Appl. 153 (2012), no. 1, 123--138.

\bibitem{Pey_Sor_2010} J. Peypouquet and S. Sorin, {\em Evolution equations for maximal monotone operators: asymptotic analysis in continuous and discrete time}, J. Convex Anal. 17 (2010), 1113--1163.

\bibitem{Rockafellar} R.-T. Rockafellar, {\em Monotone operators and the proximal point algorithm}, SIAM J. Control Optim. 14 (1976), no. 5, 877--897.

\end{thebibliography}
\end{document}